\documentstyle[11pt]{article}
 
\newcommand{\ncm}{\newcommand}
\ncm{\rncm}{\renewcommand}
\rncm{\sec}{\setc{0}\section}
\ncm{\beq}{\begin{equation}}
\ncm{\eeq}{\end{equation}}
\ncm{\bea}{\begin{eqnarray}}
\ncm{\beanon}{\begin{eqnarray*}}
\ncm{\eea}{\end{eqnarray}}
\ncm{\eeanon}{\end{eqnarray*}}
\rncm{\theequation}{\thesection.\arabic{equation}}
\ncm{\setc}[1]{\setcounter{equation}{#1}}
\newcounter{eqnr}
\newcounter{axiom}

\newenvironment{eqnarrayabc}{\stepcounter{equation}
  \setcounter{eqnr}{\value{equation}}\setc{0}
  \rncm{\theequation}{\thesection.\arabic{eqnr}\alph{equation}}
  \begin{eqnarray}}{\end{eqnarray}\setc{\value{eqnr}}}
 
\ncm{\eql}{&\quad\stepcounter{equation}(\theequation a)
           &\qquad\qquad}
\ncm{\eqr}{&\quad(\theequation b)}

\ncm{\beabc}{\begin{eqnarrayabc}}
\ncm{\eeabc}{\end{eqnarrayabc}}
\ncm{\eqboxabc}[3]{\newline\parbox[t]{1.5cm}{#1}\hfill
  \parbox[b]{12cm}{\begin{eqnarray*} #3\end{eqnarray*}}\hfill
   \parbox[b]{1.5cm}{\vspace{-0.0cm}
  \begin{eqnarrayabc}#2\end{eqnarrayabc}}\newline}
\ncm{\nn}{\nonumber \\}
 
\newtheorem{thm}{Theorem}[section]
\newtheorem{prop}[thm]{Proposition}
\newtheorem{lem}[thm]{Lemma}
\newtheorem{scho}[thm]{Scholium}
\newtheorem{defi}[thm]{Definition}
\newtheorem{coro}[thm]{Corollary}

\def\End{\mbox{End}\,}
\def\Ker{\mbox{Ker}\,}
\rncm{\Im}{\mbox{Im}}

\def\Hom{\mbox{Hom}\,}

\def\id{\mbox{id}\,}
\def\Center{\mbox{Center}\,}
\def\Hypercenter{\mbox{Hypercenter}\,}
\def\Ad{\mbox{Ad}\,}
\def\Rep{{\bf rep\,}}
\def\Mod{{\bf mod\,}}
\def\Vec{{\bf vec\,}}

\def\CA{\mbox{Center}\,A}

\def\G{{\cal G}}

\def\W{{\cal W}}

\def\I{{\cal I}}

\def\IM{\delta}    
\def\IH{I}  
\def\one{{\bf 1}}
\def\1{1}
\def\du1{{\hat\1}}
\def\o{\otimes}    
\def\oo{\times}    
\def\x{\times}     
\def\cop{\Delta}

\def\eps{\varepsilon}
\def\dueps{{\hat\eps}}
\def\duS{{\hat S}}
\def\duA{{\hat A}}
\def\tr{\mbox{tr}}
\def\c{_{(1)}}
\def\cc{_{(2)}}
\def\ccc{_{(3)}}

\def\cp{_{(1')}}
\def\ccp{_{(2')}}
\def\PL{\sqcap^L}
\def\PR{\sqcap^R}

\def\duPL{{\hat\PL}}

\def\RL{\lvec R}
\def\RR{\rvec R}
\def\RLp{\lvec{R'}}
\def\RRp{\rvec{R'}}
\def\doublehat#1{\hat{\hat{#1}}}
 
\def\lvec#1{\vbox{\ialign{##\crcr\scriptsize\leftarrowfill
            \crcr\noalign{\nointerlineskip}
            $\hfil\displaystyle{#1}\hfil$\crcr}}}
\def\rvec#1{\vbox{\ialign{##\crcr\scriptsize\rightarrowfill
            \crcr\noalign{\nointerlineskip}
            $\hfil\displaystyle{#1}\hfil$\crcr}}}
\def\lrvec#1{\vbox{\ialign{##\crcr\scriptsize\leftarrowfill
            \crcr\noalign{\vskip-12pt}\scriptsize\rightarrowfill
            \crcr\noalign{\nointerlineskip}
            $\hfil\displaystyle{#1}\hfil$\crcr}}}
\def\rlvec#1{\vbox{\ialign{##\crcr\scriptsize\rightarrowfill
            \crcr\noalign{\vskip-12pt}\scriptsize\leftarrowfill
            \crcr\noalign{\nointerlineskip}
            $\hfil\displaystyle{#1}\hfil$\crcr}}}
\def\llvec#1{\vbox{\ialign{##\crcr\scriptsize\leftarrowfill
            \crcr\noalign{\vskip-12pt}\scriptsize\leftarrowfill
            \crcr\noalign{\nointerlineskip}
            $\hfil\displaystyle{#1}\hfil$\crcr}}}
\def\rrvec#1{\vbox{\ialign{##\crcr\scriptsize\rightarrowfill
            \crcr\noalign{\vskip-12pt}\scriptsize\rightarrowfill
            \crcr\noalign{\nointerlineskip}
            $\hfil\displaystyle{#1}\hfil$\crcr}}}
 
\def\nc#1{\vbox{\ialign{##\crcr\hfil
            \crcr\noalign{\vskip-5pt}
            $\scriptstyle{#1}$\crcr}}}
\def\lc#1{\vbox{\ialign{##\crcr\tiny\leftarrowfill
            \crcr\noalign{\vskip-8pt}
            $\hfil\scriptstyle{#1}\hfil$\crcr}}}
\def\rc#1{\vbox{\ialign{##\crcr\tiny\rightarrowfill
            \crcr\noalign{\vskip-8pt}
            $\hfil\scriptstyle{#1}\hfil$\crcr}}}
\def\lrc#1{\vbox{\ialign{##\crcr\tiny\leftarrowfill
            \crcr\noalign{\vskip-12pt}\tiny\rightarrowfill
            \crcr\noalign{\vskip-8pt}
            $\hfil\scriptstyle{#1}\hfil$\crcr}}}
\def\rlc#1{\vbox{\ialign{##\crcr\tiny\rightarrowfill
            \crcr\noalign{\vskip-12pt}\tiny\leftarrowfill
            \crcr\noalign{\vskip-8pt}
            $\hfil\scriptstyle{#1}\hfil$\crcr}}}

\ncm{\fourinc}{\put(7,7){\vector(1,1){9}}
            \put(14,-13){\vector(-1,1){9}} }
\ncm{\twoinc}{\put(7,7){\vector(1,1){9}}}
\ncm{\oneinc}{\put(-5,7){\vector(-1,1){9}}}
 
\def\cros{\raise1.9pt\hbox{$\scriptscriptstyle
          > $}\!\raise1.5pt\hbox{$\scriptstyle\triangleleft\,$}}
\def\R{{I\!\!R}}
\def\C{\,{\raise1.5pt\hbox{$\scriptscriptstyle |$}
        \thinmuskip=4mu \!\!C\thinmuskip=3mu}}
\def\Z{{Z\!\!\!Z}}
\def\N{{\thinmuskip = 5.5mu I\!N\thinmuskip = 3mu}}
\def\W{{\cal W}}
\def\bra{\langle}
\def\ket{\rangle}
\def\la{\!\rightharpoonup\!}
\def\ra{\!\leftharpoonup\!}
\def\qed{\hfill {\it Q.e.d.}}
\def\Proof{{\em Proof}\,:\ }
\def\triv{D_{\varepsilon}}
\def\V{V_{\eps}}
\ncm{\Sec}{{\cal S}\!{\it ec}\,}
\ncm{\Vac}{{\cal V}\!{\it ac}\,}
\ncm{\Hyp}{{\cal H}\!{\it yp}\,}
\def\d{{\bf d}}
\def\Ind{\mbox{Index}\,}
\ncm\amalgo[1]{{\lower9pt\hbox{$\o$}\atop\raise2pt\hbox{
            $\scriptscriptstyle #1$}}}
\def\sumq{\sum_{\stackrel{q\in\Sec A}{q^L=\mu,\ q^R=\nu}}}

\begin{document}
 
\large
\title{\bf Weak Hopf Algebras II.\\
       Representation Theory, Dimensions, and the Markov Trace}
 
\author{\sc Gabriella B\"ohm $^1$,  Korn\'el Szlach\'anyi $^2$\\
\\
Research Institute for Particle and Nuclear Physics\\
H-1525 Budapest 114, P.O.B. 49, Hungary}
 
\date{}
 
\maketitle
 
\footnotetext[1]{E-mail:
  BGABR@rmki.kfki.hu\\
  Supported by the Hungarian Scientific Research Fund, OTKA --
  T 016 233}
 
\footnotetext[2]{E-mail:
  SZLACH@rmki.kfki.hu\\
  Supported by the Hungarian Scientific Research Fund,
  OTKA -- T 020 285.}

\vskip 2truecm

\begin{abstract}
If $A$ is a weak $C^*$-Hopf algebra then the category of finite
dimensional
unitary representations of $A$ is a monoidal $C^*$-category
with monoidal unit being the GNS representation $\triv$ associated
to the counit $\eps$. This category has isomorphic left dual
and right dual objects which leads, as usual, to the notion of
dimension function. However, if $\eps$ is not pure the dimension
function is matrix valued with rows and columns labelled by the
irreducibles contained in $\triv$.
This happens precisely when the inclusions $A^L\subset A$ and
$A^R\subset A$ are not connected. Still there exists a trace on
$A$ which is the Markov trace for both inclusions.
We derive two numerical
invariants for each $C^*$-WHA of trivial hypercenter. These
are the common indices $\IH$ and $\IM$, of the Haar, respectively
Markov conditional expectations of either one of the inclusions
$A^{L/R}\subset A$ and $\duA^{L/R}\subset\duA$.
In generic cases $\IH>\IM$. In the special case of weak Kac algebras
we show that $\IH=\IM$ is an integer.
\vskip 1.2truecm
\hskip 8truecm {\em Submitted to J. Algebra}
\end{abstract}

\vfill\newpage
\tableofcontents
\section{Introduction}
 
We continue the analysis of weak Hopf algebras started in
\cite{BNSz} the main issue now being the structure of
weak $C^*$-Hopf algebras. We use the notations
and terminology of \cite{BNSz} which will be referred to as I.
and the theorems, equations, etc. there will be quoted as (I.3.12)
for example.
 
Being a "quantum groupoid", i.e. a generalized concept of
symmetry,
weak Hopf algebras (WHA's) have representation categories with
monoidal product and notions of left dual and right dual
objects. In case of $C^*$-WHA's this category $\Rep A$ is a
monoidal $C^*$-category in which the left dual and right dual are
canonically isomorphic, due to the existence of the canonical
grouplike element $g$ of Prop.I.4.9. $\Rep A$ is semisimple and
the finite set $\Sec A$ of equivalence classes of irreducibles
are called the set of {\em sectors}, a term borrowed from quantum
field theory. The subset $\Vac A$ of sectors that occur in the
decomposition of the monoidal unit $\V$ of $\Rep A$ are called
{\em vacua}. This name is supported by the behaviour of general
sectors under the monoidal product: They have a groupoidlike
composition law in which $\Vac A$ plays the role of the set of
units. Thus generic sectors can be thought of as interpolating
between different vacua, we call them {\em solitons}, again by
some, however vague, quantum field theoretic motivation.
(For an approach to solitons in algebraic quantum field theory see
\cite{Fredenhagen}.)
 
As it is well known isomorphism of the left dual and right dual
allows one to introduce a faithful tracial map
$\phi_V\colon\End V\to\End\V$ for each object $V$ of $\Rep A$
which leads then to a notion of dimension $d_V$ of representations.
For uniqueness of $\phi_V$ and therefore of $d_V$ one uses a
distinguished choice of the rigidity intertwiners inherent in the
definition of duals. If the WHA $A$ is pure, i.e. its trivial
representation $\V$ is irreducible, then this choice is precisely
the {\em standard} rigidity intertwiners of \cite{Longo-Roberts}.
If $\V$ is not irreducible, i.e. decomposes into more than one
vacuum representation, then
standardness needs a modification which results in a notion of
dimension which assigns to the representation $V$ a matrix $\d_V$
the rows and columns of which are labelled by the set of vacua.
Irreducibles $q\in\Sec A$ have dimension matrices $\d_q$ which are
products of a matrix unit with a positive number $d_q$,
sometimes also called the dimension of $q$. The matrix unit
content of $\d_q$ is, however, necessary for the dimension
function $V\mapsto \d_V$ to be multiplicative and additive.
 
Sections 2 and 3 are dealing with the structure of representation
categories of WHA's, with soliton sectors, and the dimension
matrix. As a little deviation from the main course, in Subsection
3.8 we construct Frobenius-Schur
indicators for $C^*$-WHA's that has already been introduced in
\cite{FGSzV}.
 
There is an other aspect of WHA's that go well beyond their
representation categories. It is the 2-dimensional array of
inclusions one obtains from the two inclusions
$A^L\subset A\supset A^R$ by repeated applications of the Jones
construction. This is a kind of standard invariant \cite{Popa} for
a to-be-constructed depth 2 inclusion of algebras for which the
tower $A^L\subset A\subset A\cros\duA\subset\dots$ is the (first)
derived tower. In the $C^*$ setting this offers a way to describe
finite index depth 2 inclusions of von Neumann algebras (of finite
dimensional centers) as a crossed product w.r.t an action of a
$C^*$-WHA \cite{NSzW}. This is a special case of the much more
general situation considered in \cite{EV}.
 
The above mentioned array of inclusions (see Fig.2) can be thought
of as the selfintertwiner algebras of certain 1-morphisms in a
$C^*$-2-category with duals for 1-morphisms. Although this
2-category will not be made precise in this paper, it offers
a good intuitive guideline to describe the structure of WHA's
algebraically.
 
For example, we find an extension of the dimension function to
the sectors of $A^{L/R}$, which would be meaningless in the
representation category of $A^{L/R}$ since they are not
coalgebras. The dimension $\d_a$ of a sector $a\in\Sec A^L$ is
again a matrix but with rows from $\Vac A$ and columns from $\Vac
\duA$. By additivity, the dimension matrix of
$A^L\cong\oplus_aM_{n_a}$ (more precisely
of the 1-morphism the selfintertwiner algebra of which is $A^L$)
is given by $\d^L=\sum_a\,n_a\d_a$ and it plays the role of a
generator. The dimension matrix $\d^R$ of $A^R$ (i.e. that of the
1-morphism dual to that of $A^L$) is the transpose of $\d^L$ and
those of the WHA's $A$ and $\duA$ are obtained as the matrix
products
\beq
\d_A\ =\ \d^L\d^R\ ,\qquad \d_{\duA}\ =\ \d^R\d^L\ .
\eeq
The dimension matrices $\d_A$ and $\d_{\duA}$ are of course the
same that one obtains from their representation categories, as
$C^*$-WHA's.
In case of a finite group these relations become unduly trivial:
They simply say that $A$ and $\duA$ have the same dimension and
both of them are the squares of their square roots.
 
Not every triple $A^L\subset A\supset A^R$ can become
the left and right subalgebra of a $C^*$-WHA $A$. In order to
understand what restrictions this imposes on the given inclusion
triple measure theoretic concepts, such as the Haar conditional
expectations $E^{L/R}\colon A\to A^{L/R}$ and the Markov
conditional expectations $E_M^{L/R}\colon A\to A^{L/R}$, turn out
to be useful. We prove in Theorem \ref{thm Markov} that a common
Markov trace on $A$ exists for the two inclusions $A^L\subset A$
and $A^R\subset A$ implying, among others, that the inclusion
matrices of all of the connected components of $A^L\subset A$ have
the same norm. Moreover this norm and the analogue norm for $\duA$
coincide, although the inclusion matrix of $\duA^L\subset \duA$
and that of $A^L\subset A$ may be completely different.
 
The indices $\IH$ and $\IM$ of the Haar and the Markov
conditional expectations, respectively, provide "scalar" (more
precisely hypercentral) elements of $A$. They can be expressed
algebraically in terms of the integer {\em multiplicities} $n_q$
and the intrinsic {\em dimensions} $d_q$ only in special cases. We
give these special cases here:
\bea
\IH&=&\dim A^L\cdot\sum_q d_q^2\qquad\mbox{if $A$ is pure and
$S^2|_{A^L}=\id|_{A^L}$}\\
\IM&=&\sum_q\,n_qd_q\qquad\mbox{if $A$ is pure.}
\eea
For the general case see Subsections \ref{ss: Haar} and \ref{ss:
Markov}. All these formulae generalize the well known identity
$\dim A=\sum_q n_q^2$ valid for a finite group or a finite
dimensional $C^*$-Hopf algebra. The occurence of two different
indices (in general $\IH\geq\IM$) is related to non-triviality
of $S^2$, the square of the antipode. In weak Kac algebras we show
that $\IH=\IM$ and it is always an integer. In case of pure
weak Kac algebras this integer is nothing but $\dim A/\dim A^L$,
suggesting that pure weak Kac algebras might be very close to what
has been called the blowing up of (quasi-)Hopf algebras in
\cite{BSz}.

{\bf Acknowledgement}: We thank our colleague and friend, Florian
Nill for the stimulating years we have spent in writing Part I
of this paper and regret very much that he was not able to
join us in the work for Part II.

\section{Representations of WHA's}
For $A$ an associative algebra over the field $K$ let $\Mod A$
denote the category of finite dimensional left $A$-modules.
Therefore the objects of $\Mod A$ are the finite dimensional
vector spaces $V$ equipped with an action $A\ni a, V\ni v\
\mapsto\ a\cdot v\in V$ of $A$ which is nondegenerate:
$\1\cdot\ =\id_V$. Sometimes it will be convenient
to use the algebra homomorphism $D_V\colon A\to \End_K V$, the
representation on $V$, i.e. $D_V(a)v:=a\cdot v$.
The space of intertwiners (or arrows) from
the object $V$ to the object $W$ are denoted by $\Hom(V,W)$ and
consists of $K$-linear maps $T\colon V\to W$ satisfying the
intertwiner property $T(a\cdot
v)=a\cdot T(v)$, $a\in A, v\in V$. The composition of arrows
$T_1\in\Hom(V,W)$ and $T_2\in\Hom(U,V)$ are denoted by $T_1\circ
T_2$. The unit arrow at the object $V$ is $D_V(\1)$ and will
be denoted by $\one_V$.
 
In this section we will investigate the additional structure
$\Mod A$ acquires by $A$ having a weak Hopf structure.

\subsection{Monoidal structure}
The coproduct $\cop$ allows us to define a {\em
monoidal product}
of left $A$ modules and their intertwiners. At first one chooses a
strictly monoidal tensor product $\o$ in the category $\Vec K$ of
finite dimensional vector spaces over $K$. Then for two objects
$V$ and $W$ in
$\Mod A$ one makes the tensor product $V\o W$ into a left
$A$-module by setting
$a\cdot (v\o w):=x\c\cdot v\o x\cc\cdot w$. Since this module is
degenerate in general, the monoidal product $V\oo
W$ in $\Mod A$ is defined as the submodule $\cop(\1)\cdot (V\o
W)$. For intertwiners $T_i\in\Hom(V_i,W_i)$,
$i=1,2$, the monoidal product $T_1\oo T_2\in\Hom(V_1\oo V_2,W_1\oo
W_2)$ is simply the restriction of $T_1\o T_2$ onto the subspace
$V_1\oo V_2\subset V_1\o V_2$.
Coassociativity of $\cop$ and strict monoidality of $\o$
immediately imply
\bea
(T\oo R)\oo S&=&T\oo (R\oo S)\\
(T\oo R)\circ(S\oo U)&=&(T\circ S)\oo(R\circ U)\\
\one_V\oo\one_W&=&\one_{V\oo W}\ .
\eea
Although these properties are those of a strict monoidal category,
we cannot expect $(\Mod A,\,\oo)$ to be strictly monoidal
since
the monoidal unit for $\o$ (some 1-dimensional vector space) may
not belong to $\Mod A$.
A natural candidate for the {\em unit object} (or {\em monoidal
unit}) is the trivial representation $V_{\eps}$
defined in Definition I.2.13. In the sense of relaxed monoidal
categories (see \cite{MacLane}) $V_{\eps}$ is a monoidal unit if
there exist invertible arrows
$U^L_V\in\Hom(V,V_{\eps}\oo V)$, $U^R_V\in\Hom(V,V\oo V_{\eps})$,
for each $A$-module $V$ such that they are
natural in $V$,
\beq\left.\begin{array}{rcl}
(\one_{\V}\oo T)\circ U_V^L&=&U_W^L\circ T\\
(T\oo\one_{\V})\circ U_V^R&=&U_W^R\circ T
\end{array}\right\}\quad T\in\Hom(V,W)\ ,\eeq
and satisfy the triangle identities
\bea
U^L_{V}\oo \one_{W}&=&U^L_{V\oo W}\\
\one_{V}\oo U^L_{W}&=&U^R_{V}\oo \one_{W}\\
\one_{V}\oo U^R_{W}&=&U^R_{V\oo W}
\eea
for all objects $V$ and $W$.
 
\begin{prop}
The trivial left $A$-module $V_{\eps}=\duA^R$ together with
the maps
\bea
U^L_V\colon v&\mapsto&\1\c\la\du1\o\1\cc\cdot v\ \in\
V_{\eps}\oo V\\
U^R_V\colon v&\mapsto&\1\c\cdot v\o\1\cc\la\du1\ \in\
V\oo V_{\eps}
\eea
is a unit object of $(\Mod A,\oo)$.
\end{prop}
 
{\em Proof} : The arrows $U_V^L$ and $U_V^R$ are invertible arrows
with inverses
\beabc
U_V^{'L}\colon V_{\eps}\oo V\to V\ ,&\quad&
U_V^{'L}(\varphi^R\o v)=(\1\ra\varphi^R)\cdot v\\
U_V^{'R}\colon V\oo V_{\eps}\to V\ ,&\quad&
U_V^{'R}(v\o\varphi^R)=(\varphi^R\la\1)\cdot v\ ,
\eeabc
respectively. Indeed, one can easily check that $U_V^{'L}\circ
U_V^L=\one_V$, $U_V^L\circ U_V^{'L}=\one_{V_\eps\oo V}$ and similar
expressions for the right $U$-arrows. The triangle identities in
turn follow from the $\duA$-versions of axioms (A.7a--b)
and Eqns (I.2.11a--b). The details of the calculation are omitted.
For more about the monoidal structure we refer to \cite{Nill}.\qed
 
\subsection{Left duals and right duals}
In this subsection we construct left and right dual objects
in $\Mod A$ using the antipode $S$.
 
The dual space $\hat V:=\Hom_K(V,K)$ of a left $A$-module $V$ is
canonically a right $A$-module: $\bra f\cdot x,\,v\ket:=
\bra f,\,x\cdot v\ket$, $f\in\hat V,x\in
A,v\in V$. In order to make it a left $A$-module we can use either
one of the antialgebra maps $S$ or $S^{-1}$. So the {\em left
dual} module $\lvec V$ is defined to be the dual space $\hat V$
with left $A$-action $x\cdot f:= f\cdot S(x)$ and the {\em right
dual} $\rvec V$ is the same space equipped with the action $x\cdot
f:=f\cdot S^{-1}(x)$.
 
In order to establish $V\mapsto \lvec V$ and $V\mapsto \rvec V$ as
the object maps of a left duality functor and a right duality
functor, respectively, we introduce the {\em left and
right rigidity intertwiners}
\beabc         \label{def: R}
\RLp_V\colon\lvec V\o V\to V_{\eps}\ ,&\quad&
     f\o v\mapsto f(\1\c\cdot v)\1\cc\la\du1\\
\RL_V\colon V_{\eps}\to V\o\lvec V\ ,&\quad&
     \varphi^R\mapsto\sum_i\,(\1\ra\varphi^R)\cdot v_i\o f^i\\
\RRp_V\colon V\o \rvec V\to V_{\eps}\ ,&\quad&
      v\o f\mapsto f(\1\cc\cdot v)\1\c\la\du1\\
\RR_V\colon V_{\eps}\to\rvec V\o V\ , &\quad&
         \varphi^R\mapsto\sum_i\,f^i\o(\varphi^R\la\1)\cdot v_i
\eeabc
where $\{v_i\}$ is a basis in $V$ and $\{f^i\}\subset\hat V$ is
its dual basis. More precisely, rigidity intertwiners are the
appropriate restrictions of the above maps to the subspaces $\lvec
V\oo V\subset \lvec V\o V,\dots$, etc.
 
\begin{prop}
For any object $V$ in $\Mod A$ the definitions (\ref{def: R}--d)
provide intertwiners
\beanon
\RL_V\in\Hom(V_{\eps},V\oo\lvec V)\ ,&\qquad&
\RLp_V\in\Hom(\lvec V\oo V,V_{\eps})\\
\RR_V\in\Hom(V_{\eps},\rvec V\oo V)\ ,&\qquad&
\RRp_V\in\Hom(V\oo\rvec V,V_{\eps})
\eeanon
satisfying the rigidity equations
\bea
U_V^{'R}\circ(\one_V\oo\RLp_V)\circ(\RL_V\oo\one_V)\circ U_V^L&=&
                  \one_V \label{left rig 1}\\
U_{\lc{V}}^{'L}\circ(\RLp_V\oo\one_{\lc{V}})\circ
   (\one_{\lc{V}}\oo\RL_V)\circ U_{\lc{V}}^R&=&\one_{\lc{V}}
                         \label{left rig 2}\\
U_V^{'L}\circ(\RRp_V\oo\one_V)\circ(\one_V\oo\RR_V)\circ U_V^R&=&
                   \one_V \label{right rig 1}\\
U_{\rc{V}}^{'R}\circ(\one_{\rc{V}}\oo\RRp_V)\circ(\RR_V\oo
                   \one_{\rc{V}})\circ
U_{\rc{V}}^L&=&\one_{\rc{V}}\ .\label{right rig 2}
\eea
\end{prop}
\Proof : The calculation proving left rigidity is this.
 
\noindent $\RL_V$ is an intertwiner:
\beanon
\RL_V(x\la\varphi^R)&=&(\1\ra(x\la\varphi^R))\cdot v_i\o f^i=
   \bra\varphi^R,\1\c x\ket\1\cc\cdot v_i\o f^i =\\
   &=&(x\c\ra\varphi^R)S(x\cc)\cdot v_i\o f^i=
   (x\c\ra\varphi^R)\cdot v_i\o f^i\cdot S(x\cc)=\\
   &=&x\cdot \RL_V(\varphi^R)
\eeanon
$\RLp_V$ is an intertwiner:
\beanon
\RLp_V(x\cdot(f\o v))&=&f(S(x\c)\1\c x\cc\cdot v)\1\cc\la\du1=
    f(\PR(\1\c x)\cdot v)\1\cc\la\du1=\\
    &=&f(\1\c\cdot v)\PL(x\1\cc)\la\du1=f(\1\c\cdot v)(x\la
    (\1\cc\la\du1))=\\
    &=&x\cdot\RLp_V(f\o v)
\eeanon
The rigidity equation (\ref{left rig 1}):
\beanon
LHS&\colon&v\mapsto \1\c\la\du1\o\1\cc\cdot v\mapsto
     (\1\ra(\1\c\la\du1))\cdot v_i\o f^i\o\1\cc\cdot v\\
   &\mapsto&S(\1\c)\cdot v_i\o f^i(\1\cp\1\cc\cdot
     v)\1\ccp\la\du1\\
   &\mapsto&((\1\ccp\la\du1)\la\1)S(\1\c)\1\cp\1\cc\cdot v=\\
   &=&v
\eeanon
the rigidity equation (\ref{left rig 2}):
\beanon
LHS&\colon&f\mapsto\1\c\cdot f\o\1\cc\la\du1\mapsto
   \1\c\cdot f\o(\1\ra(\1\cc\la\du1))\cdot v_i\o f^i\\
   &\mapsto&f(S(\1\c)\1\cp\1\cc\cdot v_i)\1\ccp\la\du1\o f^i\\
   &\mapsto&\1\cc\la\du1\o S^{-1}(\1\c)\cdot f\\
   &\mapsto&(\1\ra(\1\cc\la\du1)S^{-1}(\1\c)\cdot f=\\
   &=& f
\eeanon
The proof of right rigidity is analogous. \qed
 
\begin{coro}        \label{coro: conjfunc}
As a consequence of the rigidity equations we have the {\em left
and right duality functors} $\Mod A\to\Mod A$ mapping
$T\in\Hom(V,W)$ into $\lvec T\in\Hom(\lvec W,\lvec V)$
and $\rvec T \in \Hom(\rvec W,\rvec V)$, respectively, where
\bea
\lvec T&:=&          \label{left dual func}
U^{'L}_{\lc{V}}\circ(\RLp_{\nc{W}}\oo\one_{\lc{V}})\circ(\one_{\lc{W}}
\oo T\oo\one_{\lc{V}})
\circ(\one_{\lc{W}}\oo\RL_{\nc{V}})\circ U^R_{\lc{W}}\\
\rvec T&:=&          \label{right dual func}
U^{'R}_{\rc{V}}\circ(\one_{\rc{V}}\oo\RRp_{\nc{W}})\circ(\one_{\rc{V}}
\oo T\oo\one_{\rc{W}})
\circ(\RR_{\nc{V}}\oo\one_{\rc{W}})\circ U^L_{\rc{W}}
\eea
They are contravariant and antimonoidal and map the
$K$-space $\Hom(V,W)$ isomorphically onto $\Hom(\lvec W,\lvec V)$
and $\Hom(\rvec W,\rvec V)$, respectively.
\end{coro}
This is a fairly standard result, so the proof is omitted.
 
It is important to remark that, in spite of the complicated form
of the rigidity intertwiners, the left dual $\lvec T$ of an
intertwiner $T$ as well as its right dual $\rvec T$, if considered
merely as $K$-linear maps $\hat W\to\hat V$, coincide with the
transpose of $T$ with respect to the canonical pairing,
\beq
\bra\lvec T(f),v\ket\ =\ \bra f,T(v)\ket\ =\ \bra\rvec T(f),v\ket
\ ,\quad f\in\hat W,\ T\in\Hom(V,W),\ v\in V\ .
\eeq
This can be checked by explicit calculation using the definitions
of $\lvec T$, $\rvec T$, and those of the intertwiners involved.
 
Similar phenomenon can be observed if one compares the natural
isomorphisms $\vartheta^L_{V,W}\colon \lvec W\oo\lvec V\to\lvec{V
\oo W}$, $\vartheta^R_{V,W}\colon\rvec W\oo\rvec V\to\rvec{V\oo
W}$ in $\Mod A$ with the natural isomorphism
$\vartheta_{V,W}\colon\hat W\o\hat V\to
\widehat{V\o W}$ in $\Vec K$. As a matter of fact the rigidity
intertwiners satisfy the following monoidality relation
\beabc
\RL_{V\oo W}&=&(\one_{\nc{V}}\oo\one_{\nc{W}}\oo\vartheta_{V,W})\circ
 (\one_{\nc{V}}\oo\RL_{\nc{W}}\oo\one_{\lc{V}})\circ(U^R_{\nc{V}}
 \oo\one_{\lc{V}})\circ\RL_{\nc{V}}\\
\RLp_{V\oo W}&=&\RLp_{\nc{W}}\circ(\one_{\lc{W}}\oo U'^L_{\nc{W}})\circ
 (\one_{\lc{W}}\oo\RLp_{\nc{V}}\oo\one_{\nc{W}})\circ(\vartheta_{V,W}^{-1}
 \oo\one_{\nc{V}}\oo\one_{\nc{W}})
\eeabc
and similar equations for the right rigidity intertwiners. Therefore
the forgetful functor $\Mod A\to \Vec K$ sends $\vartheta^L_{V,W}$
and $\vartheta^R_{V,W}$ into $\vartheta_{V,W}$.
 
It is a standard consequence of the existence of left and right
duals that there are canonical natural isomorphisms
\beabc
\iota_{\nc{V}}&:=&U^{'L}_{\lrc{V}}\circ(\RRp_{\nc{V}}\oo\one_{\lrc{V}})
      \circ(\one_{\nc{V}}\oo\RL_{\rc{V}})\circ U^R_{\nc{V}}\quad
      \in\Hom(V,\lrvec{V})\label{iota}\\
\iota'_{\nc{V}}&:=&U^{'R}_{\rlc{V}}\circ(\one_{\rlc{V}}\oo\RLp_{\nc{V}})
      \circ(\RR_{\lc{V}}\oo\one_{\nc{V}})\circ U^L_{\nc{V}}\quad
      \in\Hom(V,\rlvec{V})\label{iota'}
\eeabc
Both of these arrows, if considered only as $K$-linear maps,
coincide with the natural isomorphism $V\to\doublehat{V}$
expressing reflexivity of the objects in $\Vec K$, i.e.
$\bra\iota_{\nc{V}}(v),f\ket=\bra
f,v\ket=\bra\iota'_{\nc{V}}(v),f\ket$
for all $f\in\hat V$, $v\in V$.
 
However, in general one cannot expect to have isomorphic
intertwiners $V\to\llvec V$ in $\Mod A$. Equivalently, $\lvec V$
and $\rvec V$ may not be isomorphic as $A$-modules. In special
WHA's in which the square of the antipode is inner one can still
construct natural isomorphisms $\sigma^L_V\colon V\to\llvec V$ and
$\sigma^R_V\colon V\to\rrvec V$ but these are not canonical as
long as they cannot be expressed in terms of the basic
intertwiners $U^{L/R},\RL,\RR,\dots$, etc. We shall return to this
question in case of the $C^*$-WHA's in Subsection \ref{ss: sov}
where the
situation is different due to the existence of a $^*$-operation
allowing one to build canonical isomorphisms $\gamma_V\colon\lvec
V\to\rvec V$.
 
A further consequence of the existence of rigidity intertwiners is
Frobenius reciprocity. There are two internal Hom's in $\Mod A$:
$\Hom^L(V,W):=W\oo\lvec V$ represents the functor
$Z\mapsto\Hom(Z\oo V,W)$ and $\Hom^R(V,W):=\rvec{V}\oo W$ represents
the functor $Z\mapsto\Hom(V\oo Z,W)$.
Notice that rigidity in the sense of
\cite{Deligne-Milne}, familiar in tensor and quasitensor
categories, cannot hold in $\Rep A$ since the relation
$\Hom^L(X,Y)\oo\Hom^L(V,W)\cong\Hom^L(X\oo V,Y\oo W)$ has no chance
in the lack of a braiding.

\section{Representations of $C^*$-WHA's}
\subsection{$\Rep A$ as a bundle over $\Mod A$}
 From now on the number field $K$ is the field $\C$ of complex
numbers and the WHA $A$ is assumed to be a $C^*$-WHA. A {\em
representation} of the $C^*$-WHA $A$ is a pair $(V,\,(\,,\,)_V)$
where $V$ is a finite dimensional vector space over $\C$ carrying
a left action of $A$, i.e. an object of $\Mod A$, and $(\ ,\ )_V$
is a scalar product making $V$ a Hilbert space such that the left
action of $A$ becomes a $^*$-representation: $(u,x\cdot v)_V=
(x^*\cdot u,v)_V$ for all $u,v\in V$ and $x\in A$. The
intertwiners from $(V,\,(\,,\,)_V)$ to $(W,\,(\,,\,)_W)$ are
defined to be the intertwiners from $V$ to $W$ in $\Mod A$. The
category so obtained will be denoted by $\Rep A$.
 
The forgetful
functor $\Phi\colon\Rep A\to\Mod A$ sending $(V,\,(\,,\,)_V)$ to
$V$
is faithful and full and plays the role of a bundle projection.
In this and the next subsections we use the shorthand
notation $V_1,V_2,\dots$ for objects in the fibre
$\Phi^{-1}(\{V\})$.
Later the subscripts will be omitted and $V$ also
may stand for an object in $\Rep A$.
 
Since any $A$-module can be made a
$^*$-representation by choosing an appropriate scalar product, the
fibre over any $V$ of $\Mod A$ is non-empty. Since $\Rep A$ is a
$C^*$-category, we have a new notion of isomorphism between two
representations, the {\em unitary equivalence}. Consider an
isomorphism $T\colon V\to W$ in $\Mod A$ and choose an
object $V_1$ in the fibre over $V$. Then there is precisely
one object $W_1$ in the
fibre over $W$ such that the lift of $T$ is a unitary equivalence
$T_1\colon V_1\to W_1$. We obtain immediately that the fibers,
viewed as full subcategories, over isomorphic objects are
isomorphic. Furthermore, all objects in the same fibre are
unitarily equivalent. If we fix a $V_1$ over $V$ while
allowing $T$ to run over all automorphisms $V\to V$ then the polar
decomposition $T_1=H_1U_1\colon V_1\to V_1$ yields on the one hand
all unitaries $U_1\colon V_1\to V_1$ and on the other hand
sets up a one-to-one correspondence between
the set of objects in the fibre and positive
invertible elements $H_1$ in $\End V_1$.
 
The monoidal product $V_1\oo W_1$ of $V_1$ over $V$ and $W_1$ over
$W$ is constructed as follows. One forms the tensor product of
Hilbert spaces $V_1\o W_1$ and then defines $V_1\oo W_1$ as the
image of the projection $D_{V\oo W}(\1)$ in $V_1\o W_1$.
The monoidal product of intertwiners are defined
accordingly. In this way monoidal product becomes a bifunctor
preserving the fibres in $\Rep A$. As for the unit object in $\Rep
A$ we have to choose one particular element in the fibre over
$\V$. Although all such objects are isomorphic we would like to
choose a scalar product which is given by the already existing
data in our WHA, namely by the counit.
\begin{lem}
The monoidal unit $\V$ of $\Mod A$, i.e. the left $A$-module
$_A\hat A^R$, equipped with the scalar product
$(\varphi^R,\psi^R):=\dueps(\varphi^{R*}\psi^R)$ is a
$^*$-representation. The maps $U^L_V$, $U^R_V$, ${U'}_V^L$, and
${U'}_V^R$ of $\Mod A$ lift to isometric arrows in $\Rep A$ such
that ${U'}_{V_1}^L=U_{V_1}^{L*}$ and ${U'}_{V_1}^R=U_{V_1}^{R*}$
for all $V_1$ in the fibre of $V$ and for all objects $V$ in
$\Mod A$. These isometric arrows make the unitary representation
$\V$ a unit object of $\Rep A$ (cf. Lemma I.2.12), called the
trivial representation.
\end{lem}
\Proof $\V$ is a $^*$-representation since
\beanon
(\psi_R,x\la\varphi_R)&=&\dueps(\psi_R^*(x\la\varphi_R))=
\dueps(\duS^{-1}(\psi_R^*)(x\la\varphi_R))=\\
&=&\dueps(x\la\duS^{-1}(\psi^*_R)\varphi_R)=
\dueps(\duS^{-1}(\psi^*_R)\varphi_R\ra x)=\\
&=&\dueps((\duS^{-1}(\psi_R^*)\ra x)\varphi_R)=\dueps(
(S^{-1}(x)\la\psi_R^*)\varphi_R)=\dueps((x^*\la\psi_R)^*\varphi_R)=\\
&=&(x^*\la\psi_R,\varphi_R)\ .
\eeanon
If $V_1$ is any $^*$-representation of $A$ and $u,v\in V$ then
\beanon
(\varphi_R\o u,U_{V_1}^Lv)&=&\dueps(\varphi_R^*(\1\c\la\du1))\,
(u,\1\cc\cdot
v)=\dueps((\1\c^*\la\varphi_R)^*)\,(\1\cc^*\cdot u,v)=\\
&=&\overline{\bra\varphi_R,\1\c\ket}\,(\1\cc\cdot
u,v)=((\1\ra\varphi_R)\cdot u,v)=\\
&=&({U'}_{V_1}^L(\varphi_R\o u),\,v)
\eeanon
hence $U_{V_1}^{L*}={U'}_{V_1}^L$. Since $U_V^L$ is
a bijection with ${U'}_V^L\circ U_V^L=\one_V$, its lift
$U_{V_1}^L$ is an isometry. Similar argument shows that
$U_{V_1}^R$ is an isometry, too. The validity of the triangle
equations in $\Rep A$ follow immediately from that of $\Mod A$.
\qed
 
\subsection{Duals in $\Rep A$}
For $V_1$ a finite dimensional Hilbert space we denote by $\hat
V$ its dual linear space and by $V_1\to\hat V$, $u\mapsto \bar
u$ the antilinear map defined by $\bar u(v):=(u,v)$. Let $\overline
V_1$ denote the space $\hat V$ equipped with the scalar product
$(\bar u,\bar v):=(v,u)$. In this way the isomorphism $u\mapsto
\bar u$ becomes an antilinear isometry $V\to\overline V$. If $V_1$
carries a $^*$-representation
of the $C^*$-WHA $A$, i.e. $V_1$ is an object of $\Rep A$, then
there are two natural left $A$-module structures $\lvec{V}$
and $\rvec{V}$ on $\hat V$ (see Section 2) but neither
of them is a $^*$-representation on $\overline V_1$. If we insist on
having duality functors in $\Rep A$ that are obtained by lifting
the duality functors of $\Mod A$ then we need
to modify the scalar product on $\overline V_1$ and must not change its
$A$-module structure. So let $\lvec{V_1}$ and $\rvec{V_1}$ be the
objects in the fibre of $\lvec{V}$, resp. $\rvec{V}$, with scalar
products
\beq  \label{Gamma}
(\bar u,\bar v)_{\lc{V_1}}:=(v,\Gamma_{V_1}u)\,,\quad
(\bar u,\bar v)_{\rc{V_1}}:=(v,{\Gamma'}^{-1}_{V_1}u)\,,
\eeq
where $\Gamma_{V_1}$, $\Gamma'_{V_1}$ are positive invertible
linear transformations of $V_1$ implementing $S^2$.
Lifting the left and
right rigidity intertwiners $\RL_V,\RLp_V,\RR_V,\RRp_V$ of
(\ref{def: R}--d) to $\Rep A$ we obtain
$\RL_{V_1}\colon \V\to V_1\oo\lvec{V_1},\dots$etc. satisfying
rigidity relations of the form (\ref{left rig 1}, \ref{left rig
2}, \ref{right rig 1}, \ref{right rig 2}) but now in $\Rep A$. The
corresponding
left and right duality functors $T\mapsto\lvec{T}$ and
$T\mapsto\rvec{T}$,
can then be defined by lifting formulae (\ref{left dual
func}, \ref{right dual func}) to
$\Rep A$. As in $\Mod A$ so in $\Rep A$, the left and right duals
$\lvec{T}$ and $\rvec{T}$ of an intertwiner $T\colon V_1\to W_1$,
if considered merely as maps $\hat W\to\hat V$, both coincide
with the transposed map $\hat T$ given by $\bra\hat T\hat
w,v\ket=\bra\hat w,Tv\ket$.
 
In a $C^*$-category it is natural to require that the duality
functors be $^*$-functors, i.e. $\lvec{T^*}=(\lvec{T})^*$ and
$\rvec{T^*}=(\rvec{T})^*$. This implies strong restrictions on
the $\Gamma_{V_1}$ and $\Gamma'_{V_1}$ in (\ref{Gamma}). For the
left dual, for example,
this leads to that $T\circ \Gamma_{V_1}=\Gamma_{W_1}\circ T$
must hold for all $T\colon V_1\to W_1$. This implies two things.
On the one hand $\Gamma_{V_1}$ has to be constant on the fibre,
and on the other hand it is natural in $V$. Similar conclusions
hold for the right dual. Finally we conclude that there exist
positive invertible elements $g_l,g_r\in A$, both of them
implementing $S^2$, such that the scalar products on all
objects $V_1$ can be written as
\beq
(\bar u,\bar v)_{\lc{V_1}}:=(v,g_l\cdot u)\,,\quad
(\bar u,\bar v)_{\rc{V_1}}:=(v,g_r^{-1}\cdot u)\,.
\eeq
The elements $g_l$ and $g_r$ will be called the {\em left metric}
and the {\em right metric}, respectively.
 
Using the $*$-operation one has more canonical arrows to build out
of the $U^{L/R}$ and $\RL,\RR$ intertwiners than it was possible
in $\Mod A$. In particular the intertwiners
\beq      \label{eq: cangamma}
\gamma_{\nc{V}}:=U_{\rc{V}}^{R\,*}\circ(\one_{\rc{V}}\oo
                 {\RL_{\nc{V}}}^*)\circ
          (\RR_{\nc{V}}\oo\one_{\lc{V}})\circ U_{\lc{V}}^L\ \in\
                 \Hom(\lvec V,\rvec V)
\eeq
are the components of a natural isomorphism between the left
dual and right dual functors. Therefore the intertwiners
\bea
\sigma^L_V&:=&\lvec\gamma_V\circ\iota_{\nc{V}}\ \colon V\to\llvec V\\
\sigma^R_V&:=&\rvec{\gamma_V^*}\circ\iota'_{\nc{V}}\ \colon V\to\rrvec V
\eea
provide canonical natural isomorphisms establishing reflexivity in
$\Rep A$. More precisely, they make the dual object functors
in $\Rep A$ reflexive in the sense of $\C$-linear categories. In
case of $C^*$-categories one requires also that $\sigma_V^L$,
$\sigma^R_V$ be isometries.
 
In the next subsections we study the question how to choose the
metrics $g_l$ and $g_r$ in order for
\begin{itemize}
\item the natural isomorphisms $\gamma_V$, $\sigma^{L/R}_V$, and
$\iota_{\nc{V}}$, $\iota'_{\nc{V}}$ to be isometries,
\item $\gamma_V$ to satisfy sovereignty in the sense of
\cite{Yetter},
\item and the rigidity intertwiners (\ref{def: R}-d) to be
standard in the sense of \cite{Longo-Roberts}.
\end{itemize}
We shall see that the above unitarity, sovereignty and
standardness conditions can be satisfied by unique $g_l$ and
$g_r$ and lead to a distinguished choice of the left dual
and right dual objects in $\Rep A$.
 
\subsection{Sovereignty}    \label{ss: sov}
 
A natural equivalence of the fibre preserving $^*$-functors
$\lvec{(\ )}$ and $\rvec{(\ )}$ is a natural isomorphism
$\gamma\colon\lvec{(\ )}\to\rvec{(\ )}$ in $\Mod A$ with all
of its components lifting to isometries $\gamma_{V_1}\colon
\lvec{V_1}\to\rvec{V_1}$. The intertwiner property
$\gamma_V(x\cdot\lvec v)=x\cdot\gamma_V(\lvec v)$, $x\in A,\lvec
v\in\lvec V$, implies that $\gamma_V$ is the tranpose of a map
$\gamma_V^t\in\End_{\C}V$ implementing $S^2$, i.e.
$\gamma_V^t(x\cdot v)=S^2(x)\cdot\gamma_V^t(v)$, $x\in A,v\in V$.
Naturality in $V$, together with semisimplicity of $A$, leads to
that $\gamma_V^t(v)=g'\cdot v$ where $g'\in A$ is
independent of $V$ and implements $S^2$. Therefore
\beq               \label{eq: gam}
\gamma_V(\lvec{v})\ =\ \rvec v\cdot g'\ =\ \rvec{g'^*\cdot v}\ ,
\qquad v\in V\ .
\eeq
Here we adopted the convention that the antilinear map $v\mapsto
\bar v$ is denoted by $v\mapsto\lvec v$ if the image is considered
to be the $A$-module $\lvec V$ and by $v\mapsto\rvec v$ if the
image is $\rvec V$. Of course, neither of these maps are
$A$-module maps:
\beq
x\cdot\lvec v\ =\ \lvec{S(x)^*\cdot v}\,,\qquad
x\cdot\rvec v\ =\ \rvec{S(x^*)\cdot v}\,.
\eeq
Now the condition for $\gamma_V$ to lift to an isometry is that
\beanon
 (\gamma_{\nc{V}}(\lvec u),\gamma_{\nc{V}}(\lvec v))_{\rc{V}}&=&
 (\rvec{g'^*\cdot u},\rvec{g'^*\cdot v})_{\rc{V}}=\\
 =(g'^*\cdot v,g_r^{-1}g'^*\cdot u)_{\nc{V}}&=&
 (v,g'g_r^{-1}g'^*\cdot u)_{\nc{V}}
\eeanon
be equal to $(\lvec u,\lvec v)_{\lc{V}}=(v,g_l\cdot u)_{\nc{V}}$
for all $u,v\in V$, i.e.
\beq                 \label{eq: unitarygamma}
g'^*g'\ =\ g_lg_r\ .
\eeq
 
By definition a natural isomorphism $\gamma\colon
\lvec{(\ )}\to\rvec{(\ )}$ is {\em monoidal} if
\bea\label{monoidality of gamma}
\vartheta^R_{V,W}\circ \gamma_W\oo\gamma_V\ =\
\gamma_{V\oo W}\circ\vartheta^L_{V,W}
\eea
and {\em sovereign} \cite{Yetter} if it is monoidal and
satisfies\footnote{We have relaxed Yetter's
condition of $\gamma_{\V}$ to be an identity arrow.}
\beq  \label{eq: sovdiag}
\gamma_{\lc{V}}^{-1}\circ\iota'_{\nc{V}}\ =\ \lvec\gamma_V
\circ\iota_{\nc{V}}\ .
\eeq
Here $\iota'_V\colon V\to\rlvec{V}$ and $\iota_V\colon
V\to\lrvec{V}$ are the lifts to $\Rep A$ of the natural
isomorphisms introduced in (\ref{iota}-b). Using (\ref{eq: gam})
the monoidality condition (\ref{monoidality of gamma}) and the
sovereignty condition (\ref{eq: sovdiag})
translate respectively to the following conditions on $g'$:
\beq
\cop(g')\ =\ (g'\o g')\cop(\1)\,,\qquad S(g')\ =\ g'^{-1}\ .
\eeq
Such grouplike elements exist in any $C^*$-WHA therefore
sovereignty natural isomorphisms exist in $\Rep A$, i.e. $\Rep A$
is sovereign. It would be tempting to choose for $g'$ the
canonical grouplike element $g$ of Proposition I.4.4.
 From the point of view of standardness, however, an other choice
will be more natural.
 
Once a choice of the dual objects is made, i.e. $g_l$ and $g_r$
have been fixed, then formula (\ref{eq: cangamma}) yields a
canonical choice for $\gamma_V$. From now on $\gamma_V$ will
always denote this natural isomorphism.
 
Of course, $\gamma_V^{*-1}$ would also be equally good. So we
require $\gamma_V$ to be isometric. In order to see what this
requirement means for the $g'$ underlying expression
(\ref{eq: cangamma}) we compute its matrix elements
\beanon
&&(\rvec v,\gamma_{\nc{V}} (\lvec u))_{\rc{V}}\ =\
\left((\one_{\rc{V}}\oo\RL_{\nc{V}})\circ U_{\rc{V}}^R(\rvec
v)\ ,\ (\RR_{\nc{V}}\oo\one_{\lc{V}})\circ U_{\lc{V}}^L(
\lvec u)\right)_{\rc{V}\oo V\oo\lc{V}}=\\
&=&\left((\one_{\rc{V}}\oo\RL_{\nc{V}})(\1\c\cdot\rvec
v\o\1\cc\la\du1)\ ,\ (\RR_{\nc{V}}\oo\one_{\lc{V}})(
\1\cp\la\du1\o\1\ccp\cdot\lvec u)\right)_{\rc{V}\oo V\oo\lc{V}}=\\
&=&\sum_{i,j}\left(\rvec v\cdot S^{-1}(\1\c)\o\1\cc\cdot v_i\o
\lvec v_i\ ,\ \rvec v_j\o\1\cp\cdot v_j\o\lvec u\cdot S(\1\ccp)
\right)_{\rc{V}\oo V\oo\lc{V}}=\\
&=&\sum_i (\rvec v\cdot S^{-1}(\1\c)\,,\,\rvec
v_i\cdot\1\cp)_{\rc{V}}
\ (\lvec v_i\cdot\1\cc\,,\,\lvec u\cdot S(\1\ccp))_{\lc{V}}=\\
&=&(\lvec{\1\cc^*\1\cp S^{-1}(\1\c)^*\cdot v}\ ,\
\lvec{S(\1\ccp)^*\cdot u})_{\lc{V}}\ =\ (g_l\cdot
u\,,\,g_r^{-1}\cdot v)_{\nc{V}}=\\
&=&(\rvec v\,,\,\rvec{g_l\cdot u})_{\rc{V}}
\eeanon
This proves that expression (\ref{eq: cangamma}) corresponds to
the choice $g'=g_l$. By Eqn (\ref{eq: unitarygamma}) this is
unitary iff also $g_r=g_l$ holds.
For this reason from now on we make the choice
$g_l=g_r=g'$ where $g'$ is positive, invertible, and implements
$S^2$. In order for $\gamma$ to be sovereign we also require
$g'$ to be grouplike.
Adjoints of rigidity intertwiners take the simple form
\bea   \label{eq: R*}
\RL_V^*&=&\RRp_V\circ(\one_V\oo\gamma_V)\\
\RR_V^*&=&\RLp_V\circ(\gamma_V^{-1}\oo\one_V)\ .
\eea
Finally we remark that together with unitarity of $\gamma$ we also
have
\begin{scho}                           \label{sch: kap,sig}
The natural isomorphisms $\iota$, $\iota'$, $\sigma^L$, and
$\sigma^R$ of $\Mod A$ lift to $\Rep A$ as follows.
\bea
\iota_{\nc{V}}\colon V\to\lrvec{V}\
&&\iota'_{\nc{V}}(v)=\lrvec{g_r\cdot v}\\
\iota'_{\nc{V}}\colon V\to\rlvec{V}\
&&\iota_{\nc{V}}(v)=\rlvec{g_l^{-1}\cdot v}\\
\sigma_V^L\colon V\to\llvec{V}\
&&\sigma_V^L(v)=\llvec{v}\\
\sigma_V^R\colon V\to\rrvec{V}\
&&\sigma_V^R(v)=\rrvec{v}
\eea
They all are isometries if we set $g_l=g_r=g'$ where $g'$ is
chosen as above.
\end{scho}

\subsection{Soliton sectors}               \label{ss: soliton}
For a while we postpone the discussion of how to fix the metric
$g'$ and turn to the groupoidlike sector composition one meets in
WHA's with reducible trivial representation.
 
Let $\sum_\nu\,P_\nu$ be the decomposition of the identity arrow
of the unit object $\V$ into minimal projections in $\End\V$.
Then by Proposition I.2.15
it is a sum over vacua and $P_\nu=\triv(z^L_\nu)$, $\nu\in\Vac A$.
\begin{lem}
If $V$ is an irreducible object of $\Rep A$ then there exists one
and only one vacuum $\nu\in\Vac A$ such that $P_\nu\oo\one_V\neq
0$. This $\nu$ depends only on the equivalence class $q$ to which
$V$ belongs therefore we write $\nu=q^L$ and call it the left
vacuum of the sector $q$. Similarly, there exists one and only one
$\nu$, depending only on the class of $V$, such that $\one_V\oo
P_\nu\neq 0$. This $\nu=q^R$ is called the right vacuum of $q$.
\end{lem}
\Proof The proof for the left vacuum goes as follows. Let $V$ be
an object in $\Rep A$ then
\beq
\one_V\ =\ U_V^{L*}\circ(\one_{\V}\oo\one_V)\circ
U_V^L\ =\ \sum_{\nu\in\Vac A}\ U_V^{L*}\circ(P_\nu\oo\one_V)\circ
U_V^L
\eeq
The right hand side is a sum of mutually orthonal
projections $L_V(\nu)\in\End V$. If $V$ is irreducible
then $\End V\cong\C$ and there is a unique $\nu$ such that
$L_V(\nu)\neq 0$. For arbitrary objects $V$ and $W$ and arbitrary
$T\colon V\to W$ the naturality of the $U$
intertwiners implies that
\beq
T\circ L_V(\nu)=L_W(\nu)\circ T\ ,\qquad\nu\in\Vac A\ .
\eeq
It follows that $\nu$ is independent on the choice of the
representant $V$ within its equivalence class. \qed
 
Let us fix a set $\{V_q\}$ of representants in each class $q$ of
irreducibles. The short hand notations $\one_p,
U_p^L,\RL_p,\dots$,etc will always refer to such representants
$V_p$.
For $p,q\in\Sec A$ we consider the monoidal product $V_p\oo V_q$.
The identity
\[
\one_p\oo \one_q=\sum_{\nu\in\Vac
A}(U_p^{R*}\oo\one_q)\circ(\one_p\oo P_\nu\oo\one_q)\circ(\one_p\oo
U_q^L)
\]
tells us that $V_p\oo V_q$ is not the zero object precisely in
case of $p^R=q^L$. In particular $(\bar q)^L=q^R$ and
$(\bar q)^R=q^L$ for all $q\in\Sec A$.
If $T\colon V_r\to V_p\oo V_q$ is a non-zero intertwiner then
\beq
T\circ L_r(\nu)\ =\ T\ = L_{V_p\oo V_q}(\nu)\circ T\ =\
(L_p(\nu)\oo\one_q)\circ T
\eeq
implies that $r^L=p^L$. Similar arguments lead to that every
irreducible $r$ occuring in the product $p\oo q$ has $r^R=q^R$.
Obviously $q^L=q=q^R$ iff $q$ is a vacuum sector.
 
The irreducible sectors $q$ for which $q^L\neq q^R$ will be called
{\em soliton sectors} since they mimic the behaviour of solitons
in $1+1$-dimensional quantum field theory as long as they connect
different vacua and compose accordingly.
 
The above characterization of soliton sectors is purely
categorical. Therefore this soliton structure occurs in any
monoidal category with semisimple identity
object. For the representation category of a $C^*$-WHA there is a
simple algebraic characterization. The vacua $\mu$ are in
one-to-one correspondence with minimal projections $z_\mu^L\in
Z^L$ and also with minimal projections $z_\mu^R=S(z_\mu^L)\in
Z^R$. The left vacuum of the sector $q\in\Sec A$ is the unique
$\mu$ for which $z_\mu^L e_q=e_q$ and its right vacuum is the
unique $\nu$ for which $e_q z_\nu^R=e_q$.
Using faithfulness of $\eps|_{A^L}$ one can easily verify that
\beabc
\PL(e_q)&=&\delta_{q\in\Vac A}\sum_{p\in\Sec A,\,p^L=q}\ e_p\ =\
z_q^L\ \delta_{q\in\Vac A}      \\
\PR(e_q)&=&\delta_{q\in\Vac A}\sum_{p\in\Sec A,\,p^R=q}\ e_p\ =\
z_q^R\ \delta_{q\in\Vac A}\ .
\eeabc
 
Let $z_H$, $H\in\Hyp A$ be the minimal hypercentral projections.
$\Hyp A$ will be called the set of {\em hyperselection sectors}
of $A$. If $z_H$ is the hypercentral support of $z_\nu^L$, or,
what is the same, of $z_\nu^R$ then we write $[\nu]=H$.
As we have seen $\one_p\oo\one_q\neq 0$ implies $p^R=q^L$. Since
the left and right vacua of a sector $q$ obviously belong to the
same hyperselection sector we obtain that $V_p\oo V_q\neq 0$
occurs only for sectors $p$ and $q$ with common hypercentral
support.
 
Now assume $[\mu]=[\nu]=H$ and ask the question
whether
there exists a sector with $q^L=\mu$ and $q^R=\nu$, i.e. whether
$z_\mu^Lz_\nu^R\neq 0$. Let ${\cal V}^H$ denote the set of vacua
in $H$ and let $\mu\sim\nu$ for $\mu,\nu\in{\cal V}^H$ be the
relation that $z_\mu^Lz_\nu^R\neq 0$. Then $\mu\sim\mu$ since
$\mu\in\Sec A$ has left and right vacuum just $\mu$.
$\mu\sim\nu\Rightarrow\nu\sim\mu$ since $q$ and $\bar q$ have
left and right vacua interchanged. Finally, $\mu\sim\nu$ and
$\nu\sim\lambda$ imply $\mu\sim\lambda$ since $V_p\oo V_q\neq 0$
precisely if the middle vacua coincide. Hence $\sim$ is an
equivalence relation and one can easily check\footnote{For $C$
a class let $z_C:=\sum_{\mu\in C}z_\mu^L$. Then
$z_C=z_C \sum_{\nu\in{\cal V}^H}z_\nu^R=z_C\sum_{\nu\in C}z^R_\nu
=z_C S(z_C)$, hence $S(z_C)=z_C$.} that the sum of
$z_\mu^L$-s within an equivalence class is $S$-invariant, hence
hypercentral. So it must be the whole of ${\cal V}^H$. This proves
that the set of sectors with left vacuum $\mu$ and right vacuum
$\nu$ is non-empty precisely if $\mu$ and $\nu$ belong to the same
hyperselection sector. This will be referred to as the "fullness"
of the hypercentral blocks (see Fig. 1).

\begin{figure}
\setlength{\unitlength}{0.5mm}
\begin{picture}(120,110)(-60,0)
\put(-5,100){\vector(0,-1){100}}\put(-15,25){$q^L$}
\put(0,105){\vector(1,0){100}}\put(75,109){$q^R$}
\multiput(0,0)(20,0){6}{\line(0,1){100}}
\multiput(0,0)(0,20){6}{\line(1,0){100}}
\multiput(1,40)(60,0){2}{\line(0,1){60}}
\multiput(0,41)(0,60){2}{\line(1,0){60}}
\multiput(61,0)(40,0){2}{\line(0,1){40}}
\multiput(60,1)(0,40){2}{\line(1,0){40}}
\put(0,40.5){\line(1,0){100}}
\put(0,100.5){\line(1,0){60}}
\put(0.5,40){\line(0,1){60}}
\put(60.5,0){\line(0,1){100}}
\put(60,0.5){\line(1,0){40}}
\put(100.5,0){\line(0,1){40}}
\multiput(10,90)(20,-20){5}{\circle{3}}
\put(5,95){\circle*{3}}
\put(25,75){\circle*{3}}
\put(25,65){\circle*{3}}
\put(35,75){\circle*{3}}
\put(45,55){\circle*{3}}
\put(55,45){\circle*{3}}
\put(75,25){\circle*{3}}
\put(85,15){\circle*{3}}
\put(95,5){\circle*{3}}
\put(8,68){$\star$}
\put(28,88){$\star$}
\put(8,48){$\star$}
\put(48,88){$\star$}
\put(23,53){$\star$}\put(34,44){$\star$}
\put(43,73){$\star$}\put(54,64){$\star$}
\put(68,8){$\star$}
\put(88,28){$\star$}
\end{picture}
\caption {The sector table. The superselection sectors $q$ of $A$
are partitioned into boxes according to their left vacuum (row)
and right vacuum (column).
Vacuum sectors $\circ$, soliton sectors $\star$, and
ordinary sectors $\bullet$.
The full submatrices are the hyperselection sectors. Conjugate
pairs of sectors are found in transposed positions. The left
regular dimension matrix $d_{\mu\nu}$ can be computed as
$\sum_{q\in box}n_qd_q$ with the $box$ at the $\mu$-th row and
$\nu$-th column.}
\end{figure}
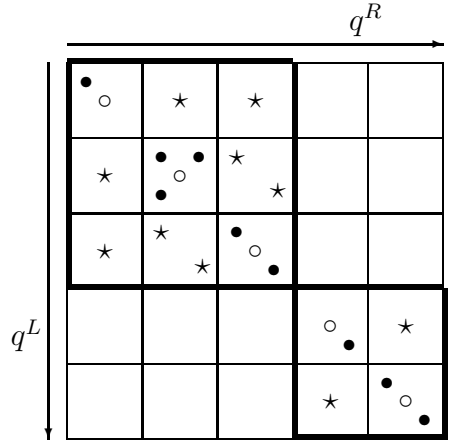

The hypercentral projections decompose our WHA into a direct sum
$\oplus_{H\in\Hyp A}\,z_HA$ of $S$-invariant
subalgebras that are subcoalgebras as well. Thus
the only "interesting" $C^*$-WHAs are those that have
trivial hypercenter. Notice, however,
that an "interesting" WHA with more than one vacuum has, in this
sense, "non-interesting" sub-WHAs for example the sum of
diagonal (i.e. non-solitonic) sectors
$A_{diag}=\oplus_{q,q^L=q^R}e_qA$. Anyhow, throughout the
paper we insist on having arbitrary hypercenter.

\subsection{Standard rigidity intertwiners}
In this subsection the selection of the
metric $g'$ will be completed by applying a further principle
called {\em standardness} \cite{Longo-Roberts}.
 
Given $g'\in A$ as in the end of Subsection \ref{ss: sov}
the left and right rigidity intertwiners (\ref{eq:
R*}) give rise to the so called left inverse and right inverse
\bea              \label{eq: phi,psi}
\phi_V\colon\ \End V\to\End\V\ ,&\quad&
\phi_V(T):=\RLp_V\circ(\gamma_V^{-1}\oo T)\circ\RR_V\ ,\\
\psi_V\colon\ \End V\to\End\V\ ,&\quad&
\psi_V(T):=\RRp_V\circ(T\oo \gamma_V)\circ\RL_V\ ,
\eea
respectively. These maps are faithful positive traces. This
can be seen either by using categorical arguments
\cite{Longo-Roberts}
or by the following direct calculation using the definitions
(\ref{def: R}--d) and (\ref{eq: gam}). At first notice that
since $\phi_V(T)=\triv(z^L)$ for some $z^L\in Z^L$,
$\phi_V(T)\varphi^R=z^L\la \varphi^R=(z^L\la\du1)\varphi^R$, so it
is sufficient to compute $\phi_V(T)$ on $\du1$.
\bea
\phi_V(T)&\colon&\du1\mapsto\sum_i\rvec{v_i}\o
v_i\ \mapsto\ \sum_i\lvec{g'^{-1}\cdot v_i}\o T(v_i)\nn
&\mapsto&\left\bra\lvec{g'^{-1}\cdot v_i},\ \1\c\cdot
T(v_i)\right\ket\ \1\cc\la\du1\ =\ \left(g'^{-1}\cdot
v_i,\ T(\1\c\cdot v_i)\right)_V\ \1\cc\la\du1\ =\nn
&=&\tr_V(T\ D_V(g'^{-1}\1\c))\ \1\cc\la\du1\nn
\phi_V(T)&=&\triv\left(\tr_V(T\
D_V(g'^{-1}\1\c))\1\cc\right) \label{eq: stand phi}
\eea
where $\tr_V$ denotes trace in the $A$-module $V$.
Similar expression can be derived for $\psi_V$. The left inverses
depend only on how we fix the freedom in $g'$. This freedom is
multiplying $g'$ with
a central positive invertible element $c$ such that $S(c)=c^{-1}$.
If $A$ is pure, i.e. $\V$ is irreducible, then such a
freedom can be eliminated by requireing $\phi_{V_q}=\psi_{V_q}$,
$q\in Sec A$, called the {\em sphericity} condition in
\cite{Westbury}. This corresponds to choosing
standard rigidity intertwiners in the sense of \cite{Longo-Roberts}
in the category $\Rep A$. If $\V$ is reducible the
sphericity condition $\phi=\psi$
cannot hold in general. Our task now is to replace this condition
with something that works for reducible identity objects as well.
 
For $V$ any object in $\Rep A$ let $T\in\End V$ and $\nu\in\Vac
A$. Then
\beq
\phi_V(U_V^{R*}\circ(T\oo P_\nu)\circ U_V^R)\ =\
P_\nu\circ\phi_V(T)
\eeq
which can be easily verified by using naturality of $U^R$ and the
triangle identities. If $V$ is irreducible, $V\cong V_q$ let us
say, then the LHS is zero for $\nu\neq q^R$ therefore $\phi_V(T)$
must be supported on $P_{q^R}$. This leads to the
\begin{defi}                    \label{def: g' and d}
Let the left and right dual objects and their rigidity
intertwiners be defined
by a common choice $g'$ for the left and right metrics. Then
$g'$ is called a {\em standard metric} if for each $q\in\Sec A$
there is a number $d_q$ such that
\beq       \label{eq: d}
\phi_{V_q}(\one_{V_q})=d_q\,P_{q^R}\ ,\qquad
\psi_{V_q}(\one_{V_q})=d_q\,P_{q^L}\ .
\eeq
The number $d_q$ is then called the {\em dimension} of the sector
$q$.
\end{defi}
 
Although our presentation is a mixture of categorical and weak
Hopf algebraic constructions, the above notion of dimension is
purely categorical. In fact Eqn (\ref{eq: d}) provides a
modification of the notion of standard left--right inverses
which works for any monoidal $C^*$-category with
duals in which the selfintertwiner space of the identity
object is finite dimensional. It is not our purpose in the present
paper to discuss standardness in general but rather to reveal the
new phenomena associated to reducibility of the identity object in
the representation categories of $C^*$-WHA's.
 
So we turn to the determination of the only possible standard metric
$g'$. Surprisingly, this $g'$ does not always coincide with the
canonical grouplike element $g$.
\begin{prop}
There exists one and only one standard metric given by the formula
\beq               \label{eq: g'}
g'\ =\ g\ k_L^{1/2}k_R^{-1/2}\ ,\ \ \mbox{where}\ \
k_L=S(k_R)=\sum_{q\in\Sec
A}\,e_q\, \dueps(\zeta_{q^L}) = \sum_{\nu\in\Vac A}\,z^L_\nu
\,\dueps(\zeta_\nu)\ .
\eeq
\end{prop}
\Proof Inserting the definitions (\ref{def: R}-d) of the the
rigidity intertwiners into (\ref{eq: phi,psi}) we obtain
\bea
\phi_V(\one_V)\equiv\RR_V^*\circ\RR_V&=&
\triv(\tr_V(g'^{-1}\1\c)\1\cc)\\
\psi_V(\one_V)\equiv\RL_V^*\circ\RL_V&=&
\triv(\tr_V(g'\1\cc)\1\c)\ .
\eea
The expressions on the RHS in the argument
of the trivial representation $\triv$ belong to $A^L$ and $A^R$,
respectively. But since $\triv$ is faithful on these subalgebras
and the left hand sides belong to $\End\V$, Proposition I.2.15
imply that
these expressions also belong to $\Center A$. If $V$ is set to
be the irreducible $V_q$ then this gives, together with the
definition of the dimension $d_q$, that
\bea
\tr_q (g'^{-1}\1\c)\1\cc&=&d_q\, z^L_{q^R} \label{eq: d_qleft}\\
\tr_q (g'\1\cc)\1\c&=&d_q\,z^R_{q^L}      \label{eq: d_qright}
\eea
where $\tr_q$ stands for $\tr_{V_q}$. Applying the counit to
these equations and utilizing the fact that $g'=gc$ with
some positive central invertible element $c=\sum_qc_qe_q$, we
obtain
\[
c_qd_q\eps(z^L_{q^R})=\tr_q(g^{-1})=\tr_q(g)=c_q^{-1}d_q\eps(z^R_{q^L})
\]
which determines $c_q$ immediately. But in order to get rid some
of the disturbing $L,R$ indices we switch to the dual WHA using
the canonical
isomorphisms of $Z^L$ and $Z^R$ with $\hat Z$ (Lemma I.2.14).
According to this Lemma there is a one-to-one
correspondence $\nu\mapsto\zeta_\nu$ of vacuum sectors of
$A$ and minimal projections of $\hat Z=\duA^L\cap\duA^R$. Hence
\beq
z^L_\nu=\1\ra\zeta_\nu\,,\quad z^R_\nu=\zeta_\nu\la\1\,,\qquad
\nu\in\Vac A\,.
\eeq
This proves formula (\ref{eq: g'}).\qed
 
Notice that by the remark after Definition I.4.11
the standard metric $g'$ is also grouplike.
\begin{coro}                 \label{cor: d}
The dimensions of the irreducible objects $V_q$ of $\Rep A$
can be expressed in terms of the weak Hopf algebraic data in the
following equivalent ways.
\beq
d_q={\tr_q(g')\over\dueps(\zeta_{q^L})}\ =\ {\tr_q(g'^{-1})\over
      \dueps(\zeta_{q^R})}
   =\tau_q/\sqrt{\dueps(\zeta_{q^L})\dueps(\zeta_{q^R})}
      \label{dq}
\eeq
where $\tau_q=\tr_q(g)$.
\end{coro}
 
In Subsection \ref{ss: dimmat} we will study the properties of
these
dimensions. Now for a little while we return to the notion of
standardness and formulate it in terms of rigidity
intertwiners with equal left dual and right dual objects that is
the common practice in $C^*$-categories. However, the content
of the next subsection is not indispensible for the rest of this
paper.

\subsection{Two-sided duals}
 
Our rigidity intertwiners $\RL,\ \RR$ were lifted from
intertwiners of $\Mod A$ where they had been associated to
different dual objects $\lvec V$ and $\rvec V$, respectively.
Although in this approach we use {\em canonical} rigidity
intertwiners built out
only of those weak Hopf algebraic data that exist for arbitrary
fields $K$,
still it is desirable to compare it with an other approach which
is more familiar in $C^*$-categories. Therefore we introduce an
"intermediate dual object" that provides a two-sided dual and find
the associated standard rigidity intertwiners.
\begin{defi}
For $V$ an object in $\Rep A$ let $\overline V$ be the dual Hilbert
space of $V$ with scalar product $(\bar u,\bar v)=(v,u)$ and left
$A$-module structure $x\cdot \bar v=\bar v\cdot
g'^{-1/2}S(x)g'^{1/2}$, where $g'$ is the standard metric.
$\overline V$ is called the {\em conjugate} of $V$.
\end{defi}
We can find isometric intertwiners
\beabc
\gamma_V^L\colon\lvec V\to\overline V\,,&&\quad \lvec v\mapsto
\bar v\cdot g'^{1/2}\\
\gamma_V^R\colon\rvec V\to\overline V\,,&&\quad \rvec v\mapsto
\bar v\cdot g'^{-1/2}
\eeabc
which satisfy $(\gamma_V^R)^{-1}\circ\gamma_V^L=\gamma_V$.
Therefore the arrows
\beabc  \label{eq: R}
R_V&:=&(\gamma_V^R\oo\one_V)\circ\RR_V\,\colon\ \V\to \overline V\oo V\\
\bar R_V&:=&(\one_V\oo\gamma_V^L)\circ\RL_V\,\colon\ \V\to
V\oo\overline V
\eeabc
satisfy the rigidity relations
\beabc
U_V^{L*}\circ(\overline R_V^*\oo\one_V)\circ(\one_V\oo R_V)
\circ U_V^R&=&\one_V\\
U_{\overline V}^{R*}\circ(\one_{\overline V}\oo \overline
R_V^*)\circ(R_V\oo\one_{\overline V})\circ U_{\overline V}^L
&=&\one_{\overline V}
\eeabc
These relations and their adjoints imply that $(\overline V,\overline
R_V,R_V^*)$ is a left dual and $(\overline V,R_V,\overline R_V^*)$ is a
right dual of $V$. This will be briefly referred to as $\overline V$ is
a {\em two-sided dual} of $V$.
 
The main advantage of the two-sided dual is that the associated
left and right dual functors coincide. This can be seen as follows.
Let $T\colon V\to W$. Then
\beanon
&&U_{\overline V}^{L*}\circ(R_W^*\oo\one_{\overline V})\circ
(\one_{\overline W}\oo
T\oo\one_{\overline V})\circ (\one_{\overline W}\oo\overline R_V)
\circ U_{\overline W}^R=\\
&=&\gamma_V^L\circ\lvec T\circ(\gamma_W^L)^{-1}=
\gamma_V^R\circ\gamma_V\circ\lvec T\circ(\gamma_V)^{-1}\circ
(\gamma_W^R)^{-1}=\gamma_V^R\circ\rvec T\circ(\gamma_W^R)^{-1}=\\
&=&U_{\overline V}^{R*}\circ(\one_{\overline V}\oo \overline
R_W^*)\circ(\one_{\overline V}\oo T\oo\one_{\overline W})\circ
(R_V\oo\one_{\overline W})\circ U_{\overline W}^L
\eeanon
We may use the notation $\overline T$ for this (left and right)
dual of $T$ and call it the {\em conjugate}. Then conjugation is a
linear
$^*$-functor, $\overline{T^*}=(\overline T)^*$. Again, as it
happened with $\lvec T$
and $\rvec T$, as a map of vector spaces, $\overline T$
coincides with the transpose of $T$ and therefore with $\lvec T$
and $\rvec T$, too. The difference is only in the $A$-module
structure and in the scalar product one puts on the vector
spaces $\hat V$ and $\hat W$.
 
The rigidity intertwiners (\ref{eq: R}-b) not only provide a
two-sided dual
but are also standard. As a matter of fact for all $q\in Sec A$
we find the normalizations
\beq
R_{V_q}^*\circ R_{V_q}=d_q\,\triv(z^L_{q^R})\,\quad
\overline R_{V_q}^*\circ \overline R_{V_q}=d_q\,\triv(z^R_{q^L})
\eeq
and for any finite direct sum
$V_i\stackrel{T_i}{\longrightarrow}
V\stackrel{T_i^*}{\longrightarrow}V_i$ of irreducibles $\{V_i\}$
\bea
R_V&=&\sum_i\ (\overline{T_i^*}\oo T_i)\circ R_{V_i}\\
\bar R_V&=&\sum_i\ (T_i\oo\overline{T_i^*})\circ\bar R_{V_i}\ .
\eea
If $\iota_V\colon V\to\hat{\hat V}$ denotes the canonical
isomorphism of finite dimensional vector spaces then it lifts to
a natural isometric isomorphism $V\to\overline{\overline V}$ in $\Rep A$.
One has the identities
\beq
R_{\overline V}=(\iota_V\oo\one_V)\circ\overline R_V\,,\quad
\overline R_{\overline V}=(\one_{\overline V}\oo\iota_V)\circ R_V\,.
\eeq
 
The left and right inverses (\ref{eq: phi,psi}) can be expressed
in terms of the two-sided duals as follows.
\bea
\phi_V(T)&=&R_V^*\circ(\one_{\bar V}\oo T)\circ R_V\\
\psi_V(T)&=&\bar R_V^*\circ(T\oo\one_{\bar V})\circ\bar R_V
\eea
Therefore they are {\em the} standard left and right inverses in
the sense of \cite{Longo-Roberts}.
 
\subsection{The dimension matrix}  \label{ss: dimmat}
 
\begin{defi}
For $V$ an object in $\Rep A$ we define its dimension matrix
$\d_V$ as follows. Its rows and columns are labelled by the
set $\Vac A$ and
\beq
\d^{\mu\nu}_V:=\sum_{\stackrel{q\in\Sec A}{q^L=\mu,\ q^R=\nu}}
\ N_V^q\,d_q
\eeq
where $N_V^q$ denotes the multiplicity of $V_q$ in $V$.
\end{defi}
For pure WHAs when $\Vac A$ has only one element this reduces to
the well known dimension formula for a reducible object $V$. The
need for introducing a matrix instead of a scalar in the non-pure
case can be seen if we ask about multiplicativity of the dimension.
Assume that a positive dimension function $V\mapsto d_V\in\R$
existed which
is multiplicative, $d_{V\oo W}=d_Vd_W$, and additive, $d_V=d_U+d_W$
for $V$ a direct sum of $U$ and $W$. Then take a soliton sector
$s$ and consider the monoidal product $V_s\oo V_s$ which is the
zero object. Hence $d_s^2=d_0=0$, a contradiction since $d_q\geq
1$ for all $q\in\Sec A$.
 
The dimension matrix can also be viewed as the set of coefficients
for the maps $Z^L\to Z^R$ and $Z^R\to Z^L$ provided by the left
inverse and the right inverse, respectively, as follows.
\bea
z_\mu^L\mapsto \phi_V(U_V^{L*}\circ(\triv(z_\mu^L)\oo\one_V)\circ
U_V^L)&=&\triv\left(\sum_\nu \d_V^{\mu\nu}z_\nu^R\right)\\
z_\mu^R\mapsto \psi_V(U_V^{R*}\circ(\one_V\oo\triv(z_\mu^R))\circ
U_V^R)&=&\triv\left(\sum_\nu z_\nu^L\d_V^{\nu\mu}\right)
\eea
The very fact that the two sets of coefficients coincide is
our standard normalization (\ref{eq: d}).
 
\begin{prop}
The dimension matrix $\d$ is an additive and multiplicative function
on the equivalence classes of objects in $\Rep A$ on which
conjugation acts by transposition. That is to say
\begin{description}
\item[i)] if $V\cong W$ then $\d_V=\d_W$,
\item[ii)] if $W$ is a direct sum of $U$ and $V$ then
$\d_W=\d_U+\d_V$,
\item[iii)] $\d_{V\oo W}=\d_V\d_W$ for all objects $V,W$,
\item[iv)] $\d_{\bar V}=\d_V^t$.
\end{description}
\end{prop}
\Proof
The only non-trivial statement is multiplicativity {\bf
(iii)} which in turn will follow from the next Lemma. \qed
\begin{lem}
The (scalar) dimension function $d_q$ given for irreducibles in
Corollary
\ref{cor: d} satisfies the following restricted multiplicativity
rule.
\beq    \label{eq: multip d}
d_p\,\delta_{p^R,q^L}\,d_q\ =\ \sum_{r\in\Sec A}\ N_{pq}^r\,d_r\ ,
\quad p,q\in\Sec A
\eeq
where $N_{pq}^r\equiv N^r_{V_p\oo V_q}$ is the multiplicity of $r$
in the product of $p$ and $q$.
\end{lem}
\Proof If $p^R\neq q^L$ then both hand sides are zero since
$N_{pq}^r=0$, $\forall r\in\Sec A$, in this case. Assume
$p^R=q^L$. Calculating $\phi_{V_p\oo V_q}(\one_p\oo\one_q)$ in two
different ways will give the required result. At first using
additivity of the left inverse yields the RHS of (\ref{eq: multip
d}). At second use grouplikeness of $g'$ to evaluate (\ref{eq:
stand phi}) with $V=V_p\oo V_q$. By means of (\ref{eq: d_qleft}) we
can write
\beanon
&&\phi_{V_p\oo V_q}(\one_p\oo \one_q)=\tr_{V_p\oo
V_q}(g'^{-1}\1\c)\triv(\1\cc)=\\
&=&\tr_p(g'^{-1}\1\cp)\tr_q(g'^{-1}\1\ccp\1\c)\triv(\1\cc)=\\
&=&d_p\,\tr_q(g'^{-1}z_{p^R}^L\1\c)\triv(\1\cc)=
d_p\,\tr_q(g'^{-1}z_{q^L}^L\1\c)\triv(\1\cc)=\\
&=&d_p\,d_q\,\triv(z_{q^R}^L)
\eeanon
where in the last equality we have utilized the fact that
$z_{q^L}^L$ is a projection containing $e_q$ as a subprojection.
\qed
 
\subsection{The Frobenius--Schur indicator}
 
In case of a finite group $G$ the Frobenius--Schur indicator is
the central element
$\iota_G:=\sum_{g\in G}g^2$ of the group algebra which takes
the values $0$ or $\pm 1/n_r$ in irreducible representations.
Non-zero values occur precisely for the selfconjugate sectors. For
$C^*$-Hopf algebras one can show that $\iota_H=h\c h\cc$, where
$h$ is the Haar measure, obeys the same properties
For a $C^*$-WHA $A$ we present here two equivalent
definitions for the Frobenius--Schur indicator, one of them is
purely categorical, the other one is Hopf algebraic. (cf.
\cite{FGSzV})
 
The
categorical definition goes as follows. Let $V$ be a selfconjugate
irreducible object in $\Rep A$. Choose an isomorphism $J\colon
V\to \overline V$. Then $\chi_V:=J^{-1}\circ\overline J\circ\iota_V\colon
V\to V$ is a number times $\one_V$ and is independent of the
choice of $J$. In particular it is isometric.
It is more tricky to show that it is selfadjoint. Using the
expression $\chi_V=U_V^{L*}\circ(\bar R_V^*\oo\one_V)\circ
(\one_V\oo J\oo J^{-1})\circ(\one_V\oo \bar R_V)\circ U_V^R$
one can verify by categorical calculus the identity $\bar \chi_V=
J\circ\chi_V^*\circ J^{-1}$. Now use the fact that $\chi_V$, being
a unitary selfintertwiner of an irreducible object, must be of the
form $u\one_V$ where $u$ is a unit length complex number.
Inserting this into our identity we obtain $u\one_{\bar V}=\bar
u\one_{\bar V}$, hence $u=\pm 1$.
This defines for each selfconjugate
sector $q=\bar q$ a number $\chi_q=\pm 1$. Extending this
definition
to $q\neq\bar q$ as $\chi_q=0$ we obtain a natural transformation
$\chi_V\colon V\to V$.
 
The Hopf algebraic definition is this. Let
$\iota_A:=h\c h\cc$. Then $\iota_A\in\Center A$ has values $0$ or
$\pm 1/\tau_r$ on irreducibles. Again the non-zero values
correspond
to selfconjugate sectors. One can show that the $\pm$ sign in
$\iota_A$ for
the sector $q$ coincides with the categorically defined $\chi_q$.
\beq
h\c h\cc\ =\ \sum_{r\in\Sec A}\ {\chi_r\over\tau_r}\,e_r
\eeq
where $\tau_r=\tr_r\,g$.
For the proof of this fact and also for clarifying
the meaning of the $\pm$ sign the following Scholium is
useful.
\begin{scho}
In a $C^*$-WHA $A$ there is a system
$\{e_r^{\alpha\beta}\,|\,r\in\Sec
A,\,\alpha,\beta=1,\dots,n_r\,\}$ of matrix units such that the
action of the antipode takes the form
\beq
g^{-1/2}S(e_r^{\alpha\beta})g^{1/2}\ =\ \left\{\begin{array}{ll}
e_{\bar r}^{\beta\alpha}&\mbox{if}\ \chi_r=0\\
e_r^{\beta\alpha}&\mbox{if}\ \chi_r=1\\
v_r e_r^{\beta\alpha} v_r^{-1}&\mbox{if}\ \chi_r=-1
\end{array}\right.
\eeq
where in the last case $n_r$ must be
even, $n_r=2k_r$, and $v_r$ in the basis
$\{e_r^{\alpha\beta}\}$ takes the form $\left(\matrix{0&I\cr
-I&0}\right)$ with $I$ denoting the $k_r\oo k_r$ unit matrix.
\end{scho}

\section{Weyl algebras as Jones extensions}
Interpreting a $C^*$-weak Hopf algebra $A$ as the algebra
generated by some set of compact, discrete coordinates and its
dual $\duA$
as the algebra generated by the associated canonical momenta
we may look for the corresponding Heisenberg or rather Weyl type
of
commutation relations. The answer is the crossed product $A\cros
\duA$ (the smash product actually) well known for finite groups
and Hopf algebras.
The novelty of the weak Hopf setting is the emergence of
an amalgamation between coordinates and momenta. We have to
identify $A^R$ with $\duA^L$ within $\W=A\cros\duA$ via the
canonical isomorphism $x^R\mapsto(\du1\ra x^R)$ of Lemma I.2.6.
 
In the Hopf algebra case the Weyl algebra $\W$ is known to
be isomorphic to
$\End_{\C}A$ which is clearly the Jones extension of the
inclusion
$\C\1\subset A$ of scalars in the Hopf algebra $A$. As a weak Hopf
generalization we will show that $\W$ is the Jones extension of
$A^L\subset A$. The non-trivial result will be that the Markov
trace of this inclusion has trace vector $t_q=f_{q^L}d_qf_{q^R}$
which is in general different from the dimension vector $d_q$.
The appearence of the positive weights $f_\nu$ on vacua $\nu$ and
the existence of a multiplicative extension of the dimension to
sectors of $A^L$ and $A^R$ reveals a genuine 2-categorical
structure underlying the $C^*$-WHA $A$. This structure enables us
to prove that each connected component of anyone of the inclusions
$A^L\subset A$, $A^R\subset A$, $\duA^L\subset\duA$, has the same
Perron-Frobenius eigenvalue, i.e. the same minimal index.
For pure $C^*$-WHA-s this minimal index takes the form
$\sum_q\,n_qd_q$ where $n_q$ is the natural number characterizing
the size of the block $q$ while $d_q\geq 1$ is the intrinsic
dimension of $q$ derived from the category $\Rep A$ in Section 3.
This generalizes the index $\sum_q d_q^2$ one obtains for
$C^*$-Hopf algebras \cite{Longo},\cite{Szymanski} in which case
$n_q=d_q$ are integers.

\subsection{The crossed product $C^*$-algebra $A\cros\duA$}
 
Any WHA $A$ is an $A^L$-$A^R$-bimodule in the obvious way. By the
canonical isomorphisms $\kappa^L_A$ and $\kappa^R_A$ of Lemma
(I.2.6) $\duA$ becomes an $A^R$-$A^L$-bimodule. One can thus form
the bimodule tensor products (or amalgamated tensor products)
$A\amalgo{A^R}\duA$ and $\duA\amalgo{A^L}A$.
\begin{defi}
The Weyl algebra $\W=\W(A)$ of a $C^*$-WHA $A$ is the crossed
product
$^*$-algebra $A\cros\duA$ with respect to the left Sweedler arrow
action of $\duA$ on $A$. This means that $\W=A\amalgo{A^R}\duA$,
as a linear space, and the multiplication and $^*$-operation are
defined respectively by
\bea
(x\o\varphi)(y\o\psi)&:=&x(\varphi\c\la y)\o\varphi\cc\psi\\
(x\o\varphi)^*&:=&\varphi^*\c\la x^*\o\varphi^*\cc\ .
\eea
\end{defi}
For showing that the above definition is independent of
the choice of the representants $x\o\varphi$ one needs only the
identities of Scholium I.2.7. It is also easy to see that $\W$
contains
$A$ and $\duA$ as unital $^*$-subalgebras. As a matter of fact
\beanon
(x\o\du1)(y\o\du1)&=&x(\du1\c\la y)\o\du1\cc=xy\c((\du1\ra
y\cc)\la\1)\o\du1=\\
&=&xy\c\PR(y\cc)\o\du1=xy\o\du1\\
(\1\o\varphi)(\1\o\psi)&=&\varphi\c\la\1\o\varphi\cc\psi=\1\o
(\du1\ra(\varphi\c\la\1))\varphi\cc\psi=\\
&=&\1\o\duPL(\varphi\c)\varphi\cc\psi=\1\o\varphi\psi
\eeanon
Identifying $x\in A$ with $x\o\du1$ and $\varphi\in\duA$ with
$\1\o\varphi$ the basic commutation relation of the Weyl algebra
reads as
\bea
\varphi x\ =\ x\c\bra x\cc,\varphi\c\ket\varphi\cc\ .
\eea
 
The following construction will show that $\W$ possesses a
faithful $^*$-representation on a Hilbert space therefore it is
actually a $C^*$-algebra.
 
The left regular $A$-module $_AA$ is a $^*$-representation if we
define the scalar product $(x,y):=\bra\hat h,x^*y\ket$ on $A$.
This Hilbert space is denoted by $L^2(A,\hat h)$.
There is an extension $\pi$ of this left regular representation to
$\W$
\beabc              \label{standard rep}
\pi(x)y&:=&xy\\
\pi(\varphi)y&:=&\varphi\la y\ .
\eeabc
This $^*$-representation is called the
{\em standard representation} of $A\cros\duA$ associated
to the Haar state. In order to prove that $\pi$ is faithful
assume that $\sum_ix_i\o\psi_i\in A\o \duA$ represents $\sum_i
x_i\psi_i\in\Ker\pi$. Then $\sum_i x_i(\psi_i\la y)=0$ for all
$y\in A$ and in particular
\beanon
&\sum_i&x_i(\psi_i\la y\cc)S^{-1}(y\c)\ =\ 0\\
&\sum_i&x_iy\cc S^{-1}(y\c)\bra\psi_i,y\ccc\ket\ =\ 0\\
&\sum_i&x_i\1\c\bra\psi_i,\1\cc y\ket\ =\ 0\\
&\sum_i&x_i\1\c\o\psi_i\ra\1\cc\ =\ 0\\
&\sum_i&x_i\1\c\o(\du1\ra\1\cc)\psi_i\ =\ 0
\eeanon
Projecting $A\o\duA$ onto $A\amalgo{A^R}\duA$ we obtain $\sum_i
x_i\psi_i=0$. Hence $\pi$ is faithful.
 
In a similar fashion one can extend the left regular
representation of $A$ to a faithful $^*$-representation $\pi'$ of
the other Weyl algebra $\duA\cros A$:
\beabc          \label{standard rep'}
\pi'(x)y&:=&xy\\
\pi'(\varphi)y&:=&y\ra \duS^{-1}(\varphi)\ .
\eeabc

\subsection{The Jones triple $A^L\subset A\subset A\cros\duA$}
 
In this subsection we show that in the faithful representation
$\pi$ the Weyl algebra is generated by $\pi(A)$ and by the
orthogonal projection $A\to A^L$. This result is a simple
application of the general method
\cite{G-H-J} in the weak Hopf environment.
Although we deviate a little bit from the standard
procedure by doing the GNS construction with respect to a
non-tracial state, namely $\hat h$, as we have discussed in the
Appendix, everything
works out as in the tracial case because the modular automorphism
of the Haar state leaves the smaller algebra (i.e. $A^L$ or $A^R$)
globally invariant.
 
Since the Haar element $\hat h$ is an idempotent we obtain for
$x^L\in A^L$ and $y\in A$ that
\beanon
(x^L,y)&=&\bra\hat h,x^{L*}y\ket=\bra\hat h,\hat h\la x^{L*}y\ket=
\bra\hat h,x^{L*}(\hat h\la y)\ket=\\
&=&(x^L,E^L(y))\ .
\eeanon
Hence the orthogonal projection onto the subspace $A^L$ is
precisely the Haar conditional expectation $E^L$ of (I.4.22).
On the other hand $E^L=\pi(\hat h)$ belongs to
$\pi(\W)$.
\begin{prop}
Let $\tilde\pi$ denote the representation of $A^{op}$ on the
Hilbert space $A$ by right multiplication, $\tilde\pi(x)y:=yx$.
Then
\beq
\pi(\W)\ =\ \tilde\pi(A^L)'\ =\ \pi(A)\,\pi(\hat h)\,\pi(A)
\eeq
Hence $\W$ is the Jones extension of $A^L\subset A$ and
the Jones projection $\pi(\hat h)$ induces the Haar conditional
expectation via
\beq                                       \label{implement}
\pi(\hat h)\pi(x)\pi(\hat h)\ =\ \pi(\hat
h)\pi(E^L(x))=\pi(E^L(x))\pi(\hat h)\ .
\eeq
\end{prop}
\Proof The identity $\varphi\la(yx^L)=(\varphi\la y)x^L$
shows that the Weyl algebra is contained in the commutant
$\tilde\pi(A^L)'$. On the one hand, the commutant is the
Jones extension which is known to be generated by $\pi(A)$ and by
the projection $e_L$ projecting onto the subspace $A^L$.
\beq
\pi(\W)\subset\tilde\pi(A^L)'=\bra\pi(A),e_L\ket
\eeq
On the other hand, we have seen above that $e_L=E^L=\pi(\hat h)$
therefore
\beq
\bra\pi(A),e_L\ket\subset\pi(\W)\,,
\eeq
which proves the main assertion. The implementation formula
(\ref{implement}) is a plain weak Hopf identity while the fact
that $\bra\pi(A),e_L\ket=\pi(A)e_L\pi(A)$ is a general result
\cite{Wenzl}. \qed
 
Analogue results hold for the right Haar conditional expectation
$E^R\colon A\to A^R$ giving rise to the Jones triple $\duA\cros
A\supset A\supset A^R$ in which the Jones extension is the other
Weyl algebra $\W(\duA)$.
 
In order to iterate the basic construction $A^L\subset A\subset
\W\subset\dots$ notice that by the left $\duA^L$-module property
of $\hat E^L$ we can apply $\id_A\o\hat E^L$ onto
$A\amalgo{A^R}\duA$ to obtain a faithful conditional expectation
$\W\to A$, also denoted by $\hat E^L$. The Jones extension
of $A\subset\W$ is the 2-fold iterated crossed product
$A\cros\duA\cros A$ in which the two copies of $A$ commute and the
first $A$ with $\duA$ satisfy $\W(A)$ commutation relations while
the last $A$ with $\duA$ satisfy $\W(\duA)$ commutation relations.
Further iterating we obtain a right growing tower of iterated
crossed products. Doing the same construction starting with
$A\supset A^R$ we obtain a left growing tower but again with the
same type of algebras. Putting the two together a 2-dimensional
array of inclusions emerges (Figure 2) in which every straight
line starting at the "sea level" is a Jones tower and the
inclusions as well as the conditional expectations commute around
each square. The algebras below the sea level are constructed by
taking intersections and their structure is governed by Lemma
I.2.14. If $A$ is a Hopf algebra then all algebras at and below
the "sea level" coincide with the number field $\C$. If $A$ is a
pure WHA such that $\duA$ is pure as well then only the $Z$
algebras below the sea are equal to $\C$. But for any WHA $A$ at
"depth 2" we reach the bottom which consists of one common copy of
Hypercenter$A$.
 
\begin{figure}
\setlength{\unitlength}{0.6mm}
\begin{picture}(200,140)(-50,-70)
\multiput(17,-60)(34,0){4}{$H$}
\multiput(0,-40)(34,0){4}{$H$}
\multiput(17,-20)(68,0){2}{$\hat Z$}
\multiput(51,-20)(68,0){2}{$Z$}
\multiput(0,0)(68,0){2}{$A^R$}
\multiput(34,0)(68,0){2}{$A^L$}
\multiput(17,20)(68,0){2}{$\hat A$}
\multiput(51,20)(68,0){2}{$A$}
\multiput(0,40)(68,0){2}{$W$}
\multiput(34,40)(68,0){2}{$\hat W$}
\multiput(17,-60)(34,0){3}{\twoinc}
\multiput(0,-40)(34,0){4}{\fourinc}
\multiput(17,-20)(34,0){3}{\fourinc}
\multiput(0,0)(34,0){4}{\fourinc}
\multiput(17,20)(34,0){3}{\fourinc}
\multiput(0,40)(34,0){4}{\fourinc}
\multiput(34,40)(34,0){3}{\oneinc}
\multiput(-27,-38)(34,0){5}{\line(1,0){24}}
\multiput(-25,3)(34,0){5}{\put(-1,0){\oval(6,6)[br]}
\multiput(5,0)(6,0){3}{\oval(6,6)[b]}
\put(23,0){\oval(6,6)[bl]} }
\put(136,0){\it sea level}\put(135,-39){\it bottom}
\put(-27,-20){\it depth}\put(-9,-1){\vector(0,-1){37}}
\put(-29,25){\it height}\put(-9,1){\vector(0,1){55}}
\put(135,-20){\sc `sea of scalars'}
\end{picture}
\caption{ The inclusion diagram of crossed product
algebras $A\cros
\hat A\cros\dots$ and $\hat A\cros A\cros\hat A\dots$ having
exactly $n$ terms ($A$ or $\hat A$) at height $n$ above sea
level. Each square is
a commuting square w.r.t. the Markov conditional expectations
$E^L_M$, $E^R_M$. Along straight lines from sea level upwards the
Jones construction is at work. Starting at depth 2 downwards the
bottom is a single copy of the hypercenter. }
\end{figure}
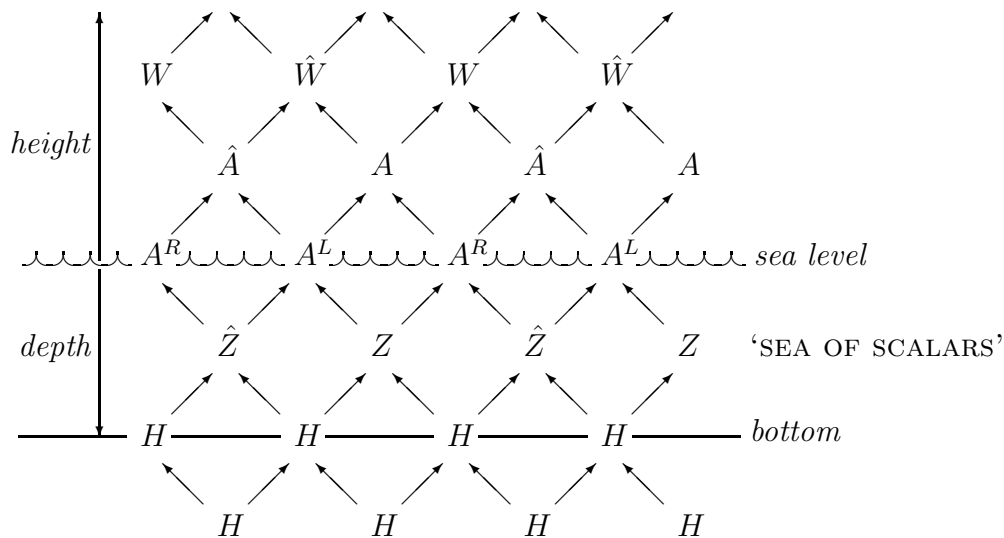
 
Exactly the above array of finite dimensional algebras of
"depth 2" arises if one considers the derived towers of a finite
index depth 2 inclusion ${\cal N}\subset{\cal M}$ of von 
Neumann algebras with relative commutant ${\cal N}'\cap{\cal
M}=A^L$ \cite{NSzW}. The members of the array can also be
interpreted
as the local algebras associated to intervals in a quantum chain.
For Hopf algebras this has been analized in \cite{NSz} and for
weak Hopf algebras in \cite{B}.

\subsection{The Haar conditional expectation}\label{ss: Haar}
 
The quasibasis of $E^L$ is by definition an element $\sum_i a_i\o
b_i\in A\o A$ satisfying $\sum_i a_i E^L(b_ix)$ $=x$ for all $x\in
A$ \cite{Watatani}. We claim that $\sum_ia_i\o
b_i=S(h\c)\o g_R^{-2}h\cc$. As a matter of fact
\beanon
S(h\c)(\hat h\la(g_R^{-2}h\cc x))&=&
S(h\c S^{-1}(x))(\hat h\la(g_R^{-2}h\cc))=\\
&=&xS(h\c)h\cc\bra\hat h,g_R^{-2}h\ccc\ket=
x\1\c\bra\hat h,g_L^{-2}h\1\cc\ket=\\
&=&x\1\c\eps(g_L^{-2}(\hat h\la h)\1\cc)\ =\ x
\eeanon
Hence we obtain for the index of $E^L$ the formula
\beq   \label{eq: ind}
\Ind E^L:=\sum_i a_ib_i=S(h\c)g_R^{-2}h\cc=\PR(g_R^{-2}h)\ \in\
A^R\cap\CA=Z^R
\eeq
which is manifestly an element of $A^R$ and belongs to $\CA$ by
the general property of the $\Ind$ \cite{Watatani}. Notice that
the connected
components of the inclusion $A^L\subset A$ are the inclusions
$z^L_\nu A^L\subset z^L_\nu A$ therefore the index being a
"scalar" would correspond to $\Ind E^L\in Z^L$. By the above
formula this is possible only if $\Ind E^L$ is hypercentral, i.e.
a true scalar in each hyperselection sector. The next
Proposition shows that this is indeed the case.
\begin{prop}
The index of the Haar conditional expectations $E^L\colon A\to
A^L$, $x\mapsto(\hat h\la x)$ and $E^R\colon A\to A^R$,
$x\mapsto(x\ra\hat h)$ is a common positive
invertible hypercentral element $\IH$ given by
\beq
\Ind E^L=\sum_{\nu\in\Vac A}z^L_\nu\ {\eps(z^L_\nu
g_L^{-2})\over \eps(z^L_\nu)}=\IH=
\sum_{\nu\in\Vac A}z^R_\nu\ {\eps(z^R_\nu
g_R^{-2})\over \eps(z^R_\nu)}=\Ind E^R\ .
\eeq
The analogue index in $\duA$ coincides with $\IH$ under the
canonical identification of the hypercenters of $A$ and $\duA$.
If the hypercenter is trivial (especially if $A$ is pure) then
the
index formula simplifies to $\IH=\1\eps(g_R^{-2})/\eps(\1)$.
\end{prop}
\Proof Using Eqns (\ref{eq: ind}) and (I.4.13) we obtain
\[\IH=\sum_q\ e_q\,\tr_q(g_L^{-1}g_R^{-1})/\tr_q(g)\]
Now we need two formulae for ratios of traces of the type
$\tr(g^\pm x^{L/R})/\tr(g)$. For that purpose
multiply Eqn (\ref{eq: d_qleft}) by $x^R\in A^R$ and Eqn (\ref{eq:
d_qright}) by $x^L\in A^L$ and then apply the counit to them. This
yields, together with (\ref{eq: g'}), the ratios
\bea                \label{tr ratio R}
{\tr_q(g^{-1}x^R)\over\tr_q(g^{-1})}&=&{\bra\zeta_{q^R},x^R\ket\over
\bra\zeta_{q^R},\1\ket}\\
{\tr_q(gx^L)\over\tr_q(g)}&=&{\bra\zeta_{q^L},x^L\ket\over
\bra\zeta_{q^L},\1\ket}\ .  \label{tr ratio L}
\eea
These formulae help to evaluate $\tr(g_L^{-1}g_R^{-1})$
in two different ways which lead immediately to the desired
expression for $\IH$. \qed
 
{\em Remark}\,: Since $E^L$ is not a trace preserving conditional
expectation (unless $\hat h$ is a trace, i.e. $S^2=\id_A$),
scalarness of $\Ind E^L$ does not necessarily imply any relation
between the connected components $z^L_\nu A^L\subset z^L_\nu A$.
(See, however, Subsection 3.4.) Thus $\IH$ may not be the minimal
index of $A^L\subset A$.
 
For later convenience we compute here the Haar conditional
expectation $E^L$ on the subalgebra $A^R$. At first notice that
$E^L(x^R)\equiv\hat h\la x^R=\1\c\bra \hat h,\,\1\cc x^R\ket$
belongs to $A^L\cap A^R=Z$. So we obtain for the Haar state on the
subalgebra $A^LA^R$ the expression
\bea
\bra\hat h,\,x^Ly^R\ket&=&\eps(x^L(\hat h\la y^R))=
\sum_{\hat\nu}\eps(x^Lz_{\hat
\nu})\frac{\eps(z_{\hat\nu}(\hat h\la
y^R))}{\eps(z_{\hat\nu})}=\nn
&=&\sum_{\hat\nu\in\Vac\duA}\,\frac{\eps(x^Lz_{\hat\nu})
\eps(z_{\hat\nu}y^R)}{\eps(z_{\hat\nu})}  \label{Haar(A^LA^R)}
\eea
and also
\beq                             \label{E^L(A^R)}
\hat h\la x^R\ =\ \sum_{\hat\nu}\, z_{\hat\nu}\,
\frac{\eps(z_{\hat\nu}x^R)}{\eps(z_{\hat\nu})}\ .
\eeq

\subsection{The Markov trace}\label{ss: Markov}
Throughout the paper $\tau$ denoted the
trace on $A$ which is related to the Haar measure by
$\tau=g_L^{-1}g_R^{-1}\la\hat h$ and has trace vector
$\tau_q=\tr_q(g)$. Since $\tau$ is faithful,
any other trace is of the form $\tau'=c\la\tau$ with $c\in\CA$.
If $\tau'$ is also faithful then the $\tau'$-preserving
conditional expectation $E^L_{\tau'}\colon A\to A^L$
can be expressed in terms of the Haar conditional expectation
as
\beq
E^L_{\tau'}(x)\ =\ E^L(rx)=E^L(xr)\,,\quad x\in A
\eeq
where the Radon-Nikodym derivative $r$ is given by (cf.
Eqn (\ref{r(s)}))
\beq
r=cg_L^{-1}g_R^{-1} E^L(cg_L^{-1}g_R^{-1})^{-1}=cg_R^{-1}(\hat
h\la cg_R^{-1})^{-1}\ . \label{q}
\eeq
The quasibasis and index of $E^L_{\tau'}$ can now be easily
obtained,
\bea
\sum_i a'_i\o b'_i&=&S(h\c)\o r^{-1}g_R^{-2}h\cc\\
\Ind E^L_{\tau'}&=&S(h\c)(\tau\la c)c^{-1}h\cc=\sum_q\,e_q\,
{\tr_q(g^{-1}c\c)\bra\tau,c\cc\ket\over\tau_q\,c_q}
                                                  \label{Ind tau'}
\eea
where $c_q$ denotes the value of $c$ in the sector $q$.
 
An important special case is the {\em standard trace} $\tau_S$
defined by the trace vector $d_q$. This is obtained by setting in
the general trace $\tau'$ the central element to be
$c=k_L^{-1/2}k_R^{-1/2}$. One would naively expect that the trace
with trace vector equal to the dimension vector is nothing else
but the Markov trace $\tau_M$, i.e. $\tau_M=\tau_S$ (up to a
hypercentral normalization). We
will see that this holds only in the absence of soliton sectors.
 
We recall that the inclusion $A^L\subset A$ is connected iff
$Z^L=\C\1$, i.e. iff $A$ is pure. Therefore in general there is no
unique Markov trace on $A$ but the Markov conditional expectation
is unique (see Definition \ref{Markov cond exp}).
Let $\Lambda=[\Lambda_{aq}]$, $a\in\Sec A^L$, $q\in \Sec A$ be the
inclusion matrix of $A^L\subset A$. Then $\Lambda^t\Lambda$
decomposes into a direct sum of irreducible matrices one for each
connected component $z_\nu^L A^L\subset z_\nu^L A$.
Hence the row (or column) indices $q$ of the matrix
$\Lambda^t\Lambda$
that belong to one and the same irreducible component can be found
in one and the same row in Figure 1. Speaking in terms of Figure
1, on the sectors $q$ of a given row there is a (up to a scalar)
unique trace vector which is the Perron-Frobenius eigenvector of
the correponding irreducible component of $\Lambda^t\Lambda$. Any
faithful trace with such a trace vector should be called a Markov
trace for the inclusion $A^L\subset A$ (or briefly a left Markov
trace) since they all share in
having the following property. The trace preserving
conditional expectation $A\to A^L$ is independent of the row by
row normalization of the trace vector and its index belongs to
$Z^L$.
It is a standard result now that the norm of this index is the
minimal
one among {\em all} conditional expectations (see Lemma A.2).
 
If we repeat this construction for the Markov trace of the
inclusion $A^R\subset A$ then the resulting trace vector will have
the freedom of an overall positive factor in each column of Figure
1. This would be a right Markov trace. Even knowing that the
inclusions $A^R\subset A$ and $A^L\subset
A$ are isomorphic via the antipode, there seems to be no reason
why the left and right Markov traces should coincide. But if they
do then they define a trace which is unique up to a scalar in each
hyperselection sector. A common left-right Markov trace would also
imply a strong relation between the disconnected parts of
$\Lambda$: They must have the same norm. Therefore that the next
Theorem is true comes as an unexpected gift of the weak Hopf
structure.
 
\begin{thm}                      \label{thm Markov}
\begin{description}
\item[i)]
There is a unique trace $\tau_M\colon A\to \C$, called {\em the
Markov trace} such that the $\tau_M$-preserving conditional
expectations $E^L_M\colon A\to A^L$ and $E^R_M\colon A\to A^R$
have equal hypercentral index $\IM$,
\beq
\Ind E^L_M\ =\ \Ind E^R_M\ =\ \IM\ \in\ \mbox{Hypercenter}\,A
\eeq
and satisfy the normalization $\tau_M(z_H)=1$ for $H\in\Hyp A$.
\item[ii)]
For any fixed hyperselection sector $H$ the
connected inclusions
\beq
z^L_\mu A^L\subset z^L_\mu A\ ,\qquad
Az^R_\mu\supset A^Rz^R_\mu,\qquad [\mu]=H
\eeq
have the same index, i.e. their inclusion matrices
$\Lambda(\mu)$ have the same
norm. This index is the value $\IM_H$ of $\IM$ on the hypersector
$H$.
\item[iii)]
$\IM$ is also equal to the norm of the dimension matrix $\d_A$
of the left regular $A$-module $_AA$. I.e. there exist numbers
$f_\mu>0$, $\mu\in\Vac A$ such that
\beq
\sum_{\nu\in\Vac A}\ d_{\mu\nu}f_\nu\ =\ \IM_{[\mu]}f_\mu
\eeq
where $d_{\mu\nu}$ stands for $\d_A^{\mu\nu}$ and $[\mu]$
denotes the hypersector of the vacuum $\mu$.
\end{description}
\end{thm}
\Proof We have seen in Subsection \ref{ss: soliton} that for
$\mu$, $\nu$
in the same hypersector $H$ there exists at least one sector $q$
with $q^L=\mu$, $q^R=\nu$. Therefore the $H$-th block of $\d_A$
is full, i.e. have strictly positive entries. In particular it is
an irreducible matrix. Let $f_\mu$, $[\mu]=H$ be a
Perron-Frobenius eigenvector and denote its eigenvalue by $\IM_H$.
 
 From the Perron-Frobenius eigenvector we can construct the central
element $f_R:=\sum_\nu\,z^R_\nu f_\nu$ $\in Z^R$ and define the
trace
$\tau_M:=c\la\tau$ with $c=c_Lc_R$, $c_R=f_R k_R^{-1/2}=S(c_L)\in
Z^R$. Then the general index formula (\ref{Ind tau'}) yields
\beanon
\Ind E^L_M&=&\sum_q\,e_q\,{\tr_q(g^{-1}\1\c)\over\tr_q(g^{-1})}
            {\bra\tau,c_R\1\cc\ket\over
            f_{q^R}\,\dueps(\zeta_{q^R})^{-1/2}}=\\
&=&\sum_q\,e_q\,{\eps(z^L_{q^R}\1\c)\over\eps(z^L_{q^R})}
{\bra\tau,c_R\1\cc\ket\over f_{q^R}\,\eps(z^L_{q^R})^{-1/2}}
=\sum_q\,e_q\,{\bra\tau,c_Rz^L_{q^R}\ket\over
f_{q^R}\,\eps(z^L_{q^R})^{1/2}}=\\
&=&\sum_\mu\ z^R_\mu\,{\bra\tau,c_Rz^L_\mu\ket\over
f_\mu\,\eps(z^L_\mu)^{1/2}}=\sum_\mu\,z^R_\mu\,\sum_{q; q^L=\mu}\
n_qd_q\,{f_{q^R}\over f_{q^L}}=\\
&=&\sum_\mu\ z^R_\mu\ \sum_\nu\,d_{\mu\nu}{f_\nu\over f_\mu}
\eeanon
which is precisely the hypercentral element $\IM$ the components
of which are the Perron-Frobenius eigenvalues $\IM_H$.
 
The restriction of $E^L_M$ onto $z^L_\mu A$ is a trace preserving
conditional expectation onto $z^L_\mu A^L$ with scalar index
$\IM_{[\mu]}$. Therefore the trace vector is the Perron-Frobenius
eigenvector and $\IM_{[\mu]}$ is the corresponding eigenvalue of
$\Lambda(\mu)^t\Lambda(\mu)$ (cf. Scholium A.4).
 
Since $\tau$ and $c$ are invariant under the antipode, so is the
Markov trace, $\tau_M\circ S=\tau_M$. Therefore
\beq  \label{S}
E^R_M\ =\ S\circ E^L_M\circ S^{-1}
\eeq
and their indices are also related by the antipode. Since the
index $\IM$ is hypercentral, they have equal index. \qed
 
The Markov trace on $A$ can now be written in the following
equivalent forms
\beq    \label{RN of tau_M/Haar}
\tau_M=f_Lf_R\la\tau_S=f_Lk_L^{-1/2}k_R^{-1/2}f_R\la\tau
      =f_Lk_L^{-1/2}g_L^{-1}g_R^{-1}k_R^{-1/2}f_R\la\hat h\ .
\eeq
Hence the trace vector of $\tau_M$ is
\beq
t_q\ :=\ f_{q^L}d_qf_{q^R}\ .
\eeq
The normalization of $\tau_M$ given in the Theorem corresponds to
the normalization of $f$ according to
\beq      \label{normalize tau_M}
\tau_M(z_H)=\sum_{q\in
H}n_qt_q=\sum_{\mu,\nu,[\nu]=H}\ d_{\mu\nu}f_\mu
f_\nu =\IM_H\sum_{\mu,[\mu]=H}\, f_\mu^2=1\ .
\eeq

\subsection{Dimensions for $A^L$, $A^R$}
The dimensions $d_q$, $q\in Sec A$ have been obtained from the
rigid monoidal structure of the category $Rep A$. Dimensions $d_a$
for the sectors $a$ of $A^L$ cannot be obtained that way since
$A^L$ is not a coalgebra hence $\Rep A^L$ is not monoidal.
However, there is an underlying 2-category ${\cal C}_A$ in "dual"
position
with respect to the representation categories in the sense that
$A$, $\duA$, $A^L,\dots$etc are selfintertwiner algebras of
certain 1-morphisms (arrows) of ${\cal C}_A$. We will not enter
into a precise construction of this 2-category here just give a
sketch of its structure on Fig. 3.
 
${\cal C}_A$ has two 0-morphisms denoted $Z$ and $\hat Z$ with
the notation refering to their selfintertwiner algebras which is
$Z=A^L\cap A^R$ and $\hat Z=\duA^L\cap \duA^R$, respectively.
There is a reducible 1-morphism $\imath$ pointing from $\hat
Z$ to
$Z$ with algebra $A^L$ and there is one, $\bar\imath$, which is
its conjugate,
pointing from $Z$ to $\hat Z$ the associated algebra of which is
$A^R$. Thus the irreducible components $a\in\Sec A^L$ have source
$a^L\in\Sec\hat Z\equiv\Vac A$ and target $a^R\in\Sec
Z\equiv\Vac\duA$. Their conjugates $b=\bar a\in\Sec A^R$ have
source $b^L\in\Vac\duA$ and target $b^R\in\Vac A$. Figure 3 is an
unfolding of this structure in order for the arrows to point to
the right and to illustrate the relation with the quantum chain of
Figure 2. The WHA $A$ corresponds to the arrow
$\imath\oo\bar\imath$
and its irreducibles to arrows $q$ connecting two vacua of $A$.
Similarly, $\duA$ is the selfintertwiner algebra of
$\bar\imath\oo\imath$ and its irreducibles $\hat q$ connect two
irreducible components of $Z$. That there are no more interesting
arrows to draw is related to the depth 2 property of Figure 2. The
graph with vertex set $\Vac A\cup\Vac\duA$ and edge set $\Sec
A^L\cup\Sec A\cup\Sec A^R\cup\Sec\duA$ will be denoted by
$\G_A$.
 
\begin{figure}
\setlength{\unitlength}{0.6mm}
\begin{picture}(200,100)(-40,25)
\multiput(0,100)(100,0){2}{\framebox(20,30)[cc]{$\hat Z$}}
\multiput(50,85)(100,0){2}{\framebox(20,30)[cc]{$Z$}}
\multiput(20,107.5)(50,0){3}{\line(1,0){30}}
\put(20,122){\line(1,0){80}}
\put(70,93){\line(1,0){80}}
\multiput(33,109)(100,0){2}{$A^L$}
\put(83,109){$A^R$}
\put(58,124){$A$}
\put(108,85){$\hat A$}
\multiput(10,30)(100,0){2}{\multiput(0,0)(0,20){3}{\circle*{3}}}
\multiput(60,40)(100,0){2}{\multiput(0,0)(0,20){2}{\circle{3}}}
\multiput(10,30)(100,0){2}{\line(5,1){48}}
\multiput(10,70)(100,0){2}{\line(5,-1){48}}
\multiput(10,70)(100,0){2}{\line(5,-3){48}}
\multiput(10,50)(100,0){2}{\line(5,1){48}}
\multiput(10,50)(100,0){2}{\line(5,-1){48}}
\put(110,30){\line(-5,1){48}}
\put(110,70){\line(-5,-1){48}}
\put(110,70){\line(-5,-3){48}}
\put(110,50){\line(-5,1){48}}
\put(110,50){\line(-5,-1){48}}
\put(9,74){$\mu$}
\put(34,68){$a$}
\put(60,32){$\hat\nu$}
\put(84,28){$b$}
\end{picture}
\caption{Schematic picture of the 1-skeleton of the 2-category
${\cal C}_A$. Upper half: The reducible 0-objects and 1-objects.
Lower half: Their irreducible content.}
\end{figure}
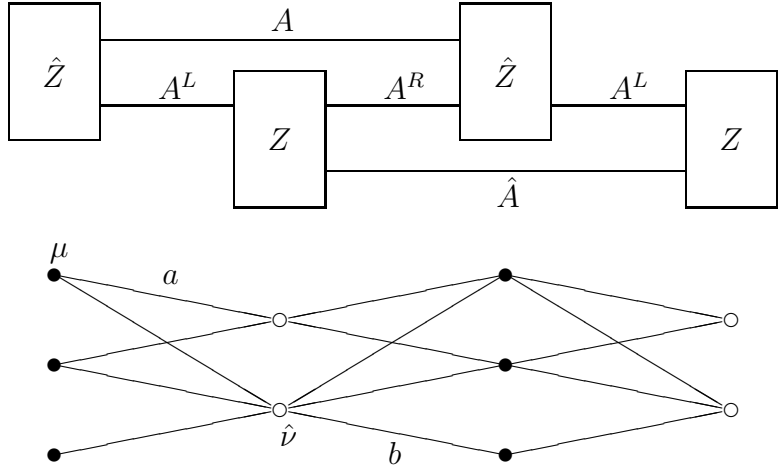
 
There are various positive functions defined on the vertices and
edges of $\G_A$.
On the vertices we have the function $k$ with values given by
the counit evaluated on the minimal projections of $Z$ and $\hat
Z$: $k_\mu=\dueps(\zeta_\mu)$, $k_{\hat\nu}=\eps(z_{\hat\nu})$.
The Perron-Frobenius eigenvector $f$ of the regular dimension
matrix (Theorem \ref{thm Markov}) determines the function $\Vac
A\ni\mu\mapsto f_\mu$ and the analogue Perron-Frobenius
eigenvector $\Vac\duA\ni\hat\nu\mapsto f_{\hat\nu}$ constructed
for $\duA$ extends $f$ to all the vertices of $\G_A$.
These functions determine the weak Hopf algebra elements
\[
\begin{array}{ll}
k=\sum_{\hat\nu}\,z_{\hat\nu}\,\eps(z_{\hat\nu})\ \in\ Z \quad&
\hat k=\sum_{\nu}\,\zeta_{\nu}\,\dueps(\zeta_{\nu})\ \in\ \hat Z\\
k_L=\1\ra\hat k\ \in\ Z^L\quad&\hat k_L=\du1\ra k\ \in\ \hat Z^L\\
k_R=\hat k\la\1\ \in\ Z^R\quad&\hat k_R=k\la\du1\ \in\ \hat Z^R\\
f=\sum_{\hat\nu}\,z_{\hat\nu}\,f_{\hat\nu}\ \in\ Z \quad&
\hat f=\sum_{\nu}\,\zeta_{\nu}\,f_{\nu}\ \in\ \hat Z \\
f_L=\1\ra\hat f\ \in\ Z^L\quad&\hat f_L=\du1\ra f\ \in\ \hat Z^L\\
f_R=\hat f\la\1\ \in\ Z^R\quad&\hat f_R=f\la\du1\ \in\ \hat Z^R
\end{array}
\]
On the edges there is the multiplicity
function $n$ : $n_a$ determines the dimension of the
irrep $a$ of $A^L$, $n_q$ that of the irrep $q$ of $A$, \dots.
For the $q$ and $\hat q$ type of edges we already have the
dimension function
\beq \label{dimensions for A, duA}
d_q=\frac{\tr_q(g^{-1})}{(k_{q^L}\,k_{q^R})^{1/2}}\ ,\qquad
d_{\hat q}=\frac{\tr_{\hat q}(\hat g^{-1})}{(k_{\hat q^L}\,
k_{\hat q^R})^{1/2}}\ .
\eeq
We now propose an extension of $d$ to the $a$ and $b$
type of edges. But before doing that let us decompose the
standard metric into left and right components as
\beq
g'\ =\ g'_L(g'_R)^{-1}\,\qquad g'_L=S(g'_R)=k_L^{1/2}g_Lk^{1/2}\ .
\eeq
Analogue formula holds for the dual standard metric $\hat g'$.
 
\begin{thm} \label{thmda}
\begin{description}
\item[i)] There exist unique positive numbers $\{d_a\}_{a\in \Sec A^L}$
$\{d_b\}_{b\in \Sec A^R}$ such that
with $N^{ab}_q$ denoting the multiplicity of the simple
algebra $e_a^LA^L\o e_b^RA^R$ in $e_qA$ we have the
multiplicativity rule
\beq            \label{d_a multip}
\sum_{q\in\Sec A}\,N_q^{ab}\,d_q\ =\ d_a\,\delta_{a^R,b^L}\,d_b
\eeq
and the conjugation rule
\beq \label{d_a conj} d_{\bar a}=d_a.\eeq
Their explicit form is given by
\beabc
d_a&=&\eps(e_a^L\,g'^{-1}_L)/n_a\ =\ \dueps(\hat e_a^R\,\hat
g'^{-1}_R)/n_a \label{da}\\
d_b&=&\eps(e_b^Rg'^{-1}_R)/n_b\ =\ \dueps(\hat e_b^L\,\hat
g'^{-1}_L)/n_b \label{db}\
\eeabc
where the minimal central idempotents of $A^L$, $A^R$,\dots etc
are related by $\hat e_a^R=e^L_a\la\du1$, $e_a^L=S(e^R_{\bar a})$,
$e^R_b=\hat e_b^L\la\1$.
Then it follows that also
\beq            \label{d_dua multip}
\sum_{{\hat q}\in\Sec \duA}\,N_{\hat q}^{ba}\,d_{\hat q}\ =
\ d_b\,\delta_{b^R,a^L}\,d_a
\eeq
holds with $N_{\hat q}^{ba}$ denoting the multiplicity of the
simple algebra
${\hat e}_b^L\duA^L\o {\hat e}_a^L\duA^R$ in $e_{\hat q}\duA$.
\item[ii)]
 Introduce the dimension matrix of $A^L$, respectively
$A^R$ by the formulae
\beabc
d_{\mu\hat\nu}&:=&\sum_{\stackrel{a\in\Sec A^L}{a^L=\mu,\
a^R=\hat\nu}} n_ad_a\ \equiv\ \eps(z^L_\mu g'^{-1}_Lz_{\hat\nu})\\
d_{\hat\mu\nu}&:=&\sum_{\stackrel{b\in\Sec A^R}{b^L=\hat\mu,\
b^R=\nu}} n_bd_b\ \equiv\ \eps(z_{\hat\mu}g'^{-1}_Rz_\nu^R)
\ =\ d_{\nu\hat\mu}.
\eeabc
Then the dimension matrices $\d_A$ and $\d_{\duA}$ of the
left regular modules of $A$, respectively $\duA$ can be expressed
as
\beq        \label{d=d^2}
d_{\mu\nu}\ =\ \sum_{\hat\rho\in\Vac\duA}\
d_{\mu\hat\rho}d_{\hat\rho\nu}\,,\quad
d_{\hat\mu\hat\nu}\ =\ \sum_{\rho\in\Vac A}\
d_{\hat\mu\rho}d_{\rho\hat\nu}\ .
\eeq
\item[iii)] The Perron-Frobenius eigenvectors $f$ and ${\hat f}$
are related by
\beq \label{frel}
\sum_\mu d_{{\hat{\nu}}\mu} f_\mu=
\delta_{[\hat\nu]}^{1/2} f_{\hat\nu}. \eeq
\end{description}
\end{thm}
\Proof Let us first prove uniqueness of $\{d_a\}$, $\{d_b\}$.
Suppose that there exists another solution $\{d'_a\}$, $\{d'_b\}$
of (\ref{d_a multip}, \ref{d_a conj}).
Then $d'_a/d_a=c_{a^R=b^L}=d_b/d'_b$ for some scalar
function $c_{\hat\nu}$ on $\Vac \duA$, so $d'_a=c_{a^R}d_a$,
$d'_{\bar a}=d_{\bar a}/c_{a^R}$. Taking into
account the conjugation rule we obtain $c_{\hat\nu}\equiv 1$.
 
In order to verify the solution (\ref{da}--b) we
compute the restriction of the Markov trace onto
$A^LA^R$ using (\ref{Haar(A^LA^R)}).
\beanon
\bra\tau_M,x^Ly^R\ket&=&\bra\hat
h,x^Lg_L^{-1}f_Lk_L^{-1/2}y^Rg_R^{-1}f_Rk_R^{-1/2}\ket=\\
&=&\sum_{\hat\nu}\eps(x^Lf_Lk_L^{-1/2}g_L^{-1}z_{\hat\nu})
\frac{1}{\eps(z_{\hat\nu})}\eps(z_{\hat\nu}g_R^{-1}k_R^{-1/2}
f_Ry^R)=\\
&=&\sum_{\hat\nu}\,\eps(x^Lf_Lg'^{-1}_Lz_{\hat\nu})\,
\eps(z_{\hat\nu}g'^{-1}_Rf_Ry^R)\ .
\eeanon
Specializing to minimal projections we obtain
\beq \label{tauMLR}
\bra\tau_M,\,e_a^L\,e_b^R\ket=f_{a^L}\,\eps(e_a^L\,g'^{-1}_L)\,
\delta_{a^R,b^L}\,\eps(g'^{-1}_Re_b^R)\,f_{b^R}\ .
\eeq
This quantity must be equal to $n_an_b\sum_q N_q^{ab}t_q$.
Taking into account the value $t_q=f_{q^L}d_qf_{q^R}$ of the trace
vector we arrive to the multiplicativity formula (\ref{d_a
multip}).
The conjugation formula is a simple consequence of the
$S$-invariance of the counit. The second expressions in
(\ref{da}--b), reflecting the symmetric roles of $A$ and
$\duA$, imply (\ref{d_dua multip}).
 
(\ref{d=d^2}) follows using (\ref{d_a multip}), (\ref{d_dua multip})
and the dimension counting formulae
\beq
\sum_{a\in\Sec A^L}\ \sum_{b\in\Sec A^R }\ N^{ab}_q\,n_an_b\ =\ n_q
\qquad
\sum_{a\in\Sec A^L}\ \sum_{b\in\Sec A^R }\ N^{ba}_{\hat q}\,n_a n_b\
=\ n_{\hat q}.
\eeq
The proof of (\ref{frel}) is now straightforward. \qed
 
The importance of the Markov index $\IM$ in $C^*$-WHA's can be
illustrated by the three-fold role in which it appears as a
Perron-Frobenius eigenvalue:
\begin{itemize}
\item $\IM$ is the PF-eigenvalue of the left regular dimension
matrices $\d_A$ and $\d_{\duA}$ (cf. Theorem \ref{thm Markov}.iii)
$$\sum_{\nu\in\Vac A} d_{\mu\nu} f_{\nu} =
\delta_{[\mu]}f_{\mu}\qquad
\sum_{\hat{\nu}\in\Vac \duA} d_{\hat{\mu}\hat{\nu}} f_{\hat{\nu}} =
\delta_{[\hat{\mu}]}f_{\hat{\mu}}\ .$$
\item $\IM$ is the PF-eigenvalue of $\Lambda^t\Lambda$
where $\{\Lambda_{aq}\}$ is the inclusion matrix of
$A^L\subset A$ (cf. Theorem \ref{thm Markov}.ii),
$$\sum_{q\in \Sec A} \Lambda_{aq}t_q=\IM_{[a]}^{1/2}t_a\qquad
\sum_{a\in \Sec A^L} t_a \Lambda_{aq} =
\IM_{[q]}^{1/2}t_q\ .$$
Here $t_a=f_{a^L}d_af_{a^R}$, so $\IM_{[a]}^{1/2}t_a$ is the trace
vector of the restriction of $\tau_M$ onto $A^L$.
\item $\IM$ is also the PF-eigenvalue of the matrix
$N_A=[N_{Aq}^r]$ where $N_{Aq}^r=\sum_p
n_pN_{pq}^r$ denotes the
multiplicity of $r$ in the monoidal product of the left regular
module $_AA$ with the sector $q$,
$$\sum_{r\in \Sec A} N_{Aq}^r t_r = \IM_{[q]} t_q\ .$$
\end{itemize}
 
The identity $S(\1\c)\1\cc=\1$ occured many times in
this paper but until now nothing
has been said about the element $\1\cc S(\1\c)\in A^L$.
After having introduced $g'_L$ and $d_a$ we are in the
position to do this.
\begin{lem} \label{RNder}
In any WHA $A$ over a field $K$ the element
$\1\cc S(\1\c)$ is the Radon-Nikodym derivative of the
left regular trace $\tr^L$ of the subalgebra $A^L$ with respect to
the nondegenerate functional $\eps|_{A^L}$, i.e.
\beq
\tr^L(x^L)=\varepsilon(x^L \1\cc S(\1\c))\qquad x^L\in A^L\ .
\eeq
If $A$ is a $C^*$-WHA then this Radon-Nikodym derivative is
positive and invertible and can be expressed as
\beq  \label{arel}
\1\cc S(\1\c)= g^{\prime -1}_L \sum_{a\in \Sec A^L} e_L^a
\frac{n_a}{d_a}\ .
\eeq
\end{lem}
\Proof Any representation $\cop(\1)=\sum_j e_j\o e^j$ with
elements $e_j\in A^R$, $e^j\in A^L$ determines a pair of dual
bases for the nondegenerate bilinear form
$(x^R,x^L)\mapsto\eps(x^Rx^L)$ on $A^R\x A^L$ of Lemma I.2.2.
Therefore
\beq
\tr^L(x^L)=\eps(\1\c x^L \1\cc)=\eps(x^L\1\cc S(\1\c))
\eeq
due to (I.2.10) and (I.2.2a).
In order to prove (\ref{arel}) notice that $x^L\1\cc S(\1\c)=
\1\cc S(\1\c)S^2(x^L)$, $x^L\in A^L$, therefore $w:=g'_L\1\cc
S(\1\c)$ belongs to $\Center A^L$. Then, denoting $\tr_a(x^L):=
{1\over n_a}\tr^L(e^L_ax^L)$,
\beq
d_a\equiv\frac{1}{n_a}\eps(e^L_a g^{\prime -1}_L)=
\frac{1}{n_a}\eps(e^L_a w^{-1}\1\cc S(\1\c))=\tr_a(w^{-1}).
\eeq
which proves (\ref{arel}). \qed
 
Applying $\tr_a$ to (\ref{arel}) we obtain
\beq
\label{da'} d_a=\frac{\tr_a(g_L^{-1})}{(k_{a^L}
k_{a^R})^{1/2}}
\eeq
a formula reminiscent to (\ref{dimensions for A, duA}).
 
{\em Remark}\,: We summarize without proof some results on the
inclusion $A^LA^R\subset A$. Recall that $A^LA^R$ is a
$C^*$-WHA by restricting the structure maps of $A$ and its
hypercentral blocks are precisely of the type discussed in
Subsection 5.2. The Haar state on $A^LA^R$ is the restriction
of $\hat h$ and the Haar index is $\Ind E^L|_{A^LA^R}=k$, i.e. the
element defined in Subsection 4.5. The map
\beq E\colon A\to A^LA^R\,,\quad E(x)=k\1\c E^L(x S^{-1}(\1\cc))
\eeq
is a conditional expectation which is the $\hat h$-preserving, the
$\tau$-preserving, the $\tau_S$-preserving, and the
$\tau_M$-preserving conditional expectation at the same time.
Its index is hypercentral, thus $E$ is {\em the} Markov
conditional expectation onto $A^LA^R$.
\bea
\Ind E&=&\sum_{\mu\in\Vac A} z_\mu^L\,\eps(z_\mu^L {g'_L}^{-2})
\ \in\ \Hypercenter A\nn
&=&\sum_{\mu\in\Vac A} z_\mu^L\cdot\sum_{a;a^L=\mu}\,d_a^2
\eea
This implies that the above sum of $d_a^2$-s is independent of
$\mu$ within a hyperselection sector $H$ and gives the square of
the norm of the inclusion matrix of $z_HA^LA^R\subset z_HA$. (As a
comparison, the Haar index $\IH$ for $A^L\subset A$ can be written
as $\sum_{a;a^L=\mu} d_a^2k_{a^R}$ and the Markov index $\IM$ is
not algebraically expressible in terms of the dimensions, either.)
Since $\tau_S$ is a Markov trace for $A^LA^R\subset A$ and has
trace vector $d_q$,
the dimension multiplicativity formula (\ref{d_a multip}) has as a
counterpart
\beq
\sum_a\sum_b\ N_q^{ab}d_ad_b\ =\ d_q\,(\Ind E)_{[q]}\ .
\eeq
 
\subsection{Temperley-Lieb projections}
 
Since the Weyl algebra $A\cros\duA$ is the common Jones extension
of the two inclusions $A^L\subset A$ and $\duA\supset\duA^R$, by
standard results \cite{G-H-J} there exist unique projections $\hat
e$ and $e$ in $A\cros\duA$ that implement the Markov conditional
expectations $E^L_M$ and $\hat E^R_M$, respectively, in the sense
of the formulae
\beq \label{TLJ1}\left.
\begin{array}{rcl}
\hat e x\hat e&=&E^L_M(x)\hat e\qquad x\in A\\
e\varphi e&=&\hat E^R_M(\varphi)e\qquad\varphi\in\duA
\end{array}\right\}\ \mbox{ within }\ A\cros \duA\ .
\eeq
The peculiarity of the smash product extension is that these Jones
projections not only belong to $A\cros\duA$ but $e\in A$ and $\hat
e\in\duA$, as well. Furthermore, they satisfy the Temperley-Lieb
relations
\beabc                 \label{TLJ2}
e\hat e e&=&\IM^{-1}\,e\\
\hat e e\hat e&=&\IM^{-1}\,\hat e
\eeabc
and, as a consequence of the manifest selfduality of these
relations, they also provide us with the Jones
projections for the Markov conditional expectations $E^R_M\colon
A\to A^R$ and $\hat E^L_M\colon\duA\to\duA^L$ the common Jones
extension of which is the other Weyl algebra $\duA\cros A$.
Therefore
\beq  \label{TLJ3} \left.
\begin{array}{rcl}
e\varphi e&=&\hat E^L_M(\varphi)e\\
\hat e x\hat e&=&E^R_M(x)\hat e
\end{array}\right\}\ \mbox{ within }\ \duA\cros A\ .
\eeq
Before proving these statements we recall that in finite
dimensional $C^*$-Hopf algebras it is well-known that $e=h$ and
$\hat e=\hat h$ are the Haar integrals of $A$ and $\duA$,
respectively.
Not too much surprisingly this is
not true in case of $C^*$-WHA's. As a matter of fact within
$A\cros\duA$ we have
\beq \label{hrel}
\hat h x\hat h= (\hat h\c\la x)\hat h\cc \hat h=
(\du1\c\hat h\la x)\du1\cc\hat h= (\du1\c\la E^L(x))\du1\cc\hat h=
E^L(x)\hat h
\eeq
hence $\hat h h\hat h = g_L^2 \hat h=\hat h g_L^2$.
Similarly, in $A\cros\duA$ we can write
\beq
h\varphi h= h\hat E^R(\varphi)
\eeq
hence $h\hat h h =h\hat g_R^2 = \hat g_R^2 h$.
 
In the next Theorem we use the notions of {\em standard
representation} $\pi_M$ of $A\cros\duA$ and {\em standard
representation} $\pi'_M$ of $\duA\cros A$ associated to the Markov
trace. Both of these representations act on the GNS Hilbert space
$L^2(A,\tau_M)$ associated to the functional $\tau_M$. They are
equivalent to the standard representations (\ref{standard rep}-b)
and (\ref{standard rep'}-b), respectively, by means of the
isometry
$U\colon L^2(A,\tau_M)\to L^2(A,\hat h)$, $x\mapsto xs^{1/2}$,
where $s=f_Lk_L^{-1/2}g_L^{-1}g_R^{-1}k_R^{-1/2}f_R$ is the
Radon-Nykodim derivative of $\tau_M$ with respect to $\hat h$
(\ref{RN of tau_M/Haar}). That is to say
$\pi_M=\Ad_{U^{-1}}\circ\pi$ and $\pi'_M=\Ad_{U^{-1}}\circ\pi'$.
 
\begin{thm}
The Radon-Nikodym derivatives of the Markov conditional
expectations with respect to the Haar ones are given by the
following formulae.
\beabc
E^L_M(x)=E^L(r_Rx),\ x\in A\,,\quad
r_R=\IM^{-1/2}\,f^{-1}k^{1/2}g_R^{-1}k_R^{-1/2}f_R\\
E^R_M(x)=E^R(r_Lx),\ x\in A\,,\quad
r_L=\IM^{-1/2}\,f_Lk_L^{-1/2}g_L^{-1}k^{1/2}f^{-1}
\eeabc
In terms of the quantities
\beabc
q_L&=&\IM^{-1/2}\,f_L^{-1}k_L^{1/2}g_L^{-1}k^{-1/2}f\ =\
\1\ra\hat r_R\\
q_R&=&\IM^{-1/2}\,fk^{-1/2}g_R^{-1}k_R^{1/2}f_R^{-1}\ =\
S(q_L)=\hat r_L\la\1\\
\hat q_L&=&\IM^{-1/2}\,\hat
f_L^{-1}\hat k_L^{1/2}\hat g_L^{-1}\hat k^{-1/2}\hat f\ =\
\du1\ra r_R\\
\hat q_R&=&\IM^{-1/2}\,\hat
f\hat k^{-1/2}\hat g_R^{-1}\hat k_R^{1/2}\hat f_R^{-1}\ =\
\hat S(\hat q_L)=r_L\la\du1
\eeabc
we define the projections
\bea
e:=q_L^{1/2}hq_L^{1/2}\,,\qquad
\hat e:=\hat q_L^{1/2}\hat h\hat q_L^{1/2}
\eea
that are the Jones projections associated to the Markov
conditional expectations in the following sense. The standard
representations $\pi_M$ of $A\cros\duA$ and $\pi'_M$ of
$\duA\cros A$ on the Hilbert space $L^2(A,\tau_M)$
send $\hat e$ to the orthogonal projection onto
the subspace $A^L$ and $A^R$, respectively:
\bea
\pi_M(\hat e)=E^L_M\ ,\qquad\pi'_M(\hat e)=E^R_M\ .
\eea
$e$ does the same after interchanging the roles of $A$ and $\duA$.
Furthermore $e$ and $\hat e$ satisfy the relations (\ref{TLJ1}),
(\ref{TLJ2}-b), and (\ref{TLJ3}).
\end{thm}
\Proof Using (\ref{r(s)}) we have
$r_R=sE^L(s)^{-1}=f_Rk_R^{-1/2}g_R^{-1}(\hat h\la
f_Rk_R^{-1/2}g_R^{-1})^{-1}$ so we need formula (\ref{E^L(A^R)}):
\bea {\hat h}\la f_Rk_R^{-1/2}g_R^{-1} &=&
\sum_{{\hat\nu}\in \Vac \duA} z_{\hat \nu}
\frac{\eps(z_{\hat\nu}g_R^{-1}k_R^{-1/2}f_R)}{k_{\hat\nu}} =
\sum_{{\hat\nu}\in\Vac\duA}\sum_{\mu\in\Vac A}z_{\hat\nu}
\frac{\eps(z_{\hat\nu}g_R^{-1}z_{\mu}^R)}{k_{\hat\nu}k_\mu^{1/2}}
f_{\mu}=\nn
&=&\sum_{\hat\nu\in\Vac\duA}\sum_{\mu\in\Vac A}\ z_{\hat\nu}
\frac{d_{\hat\nu\mu}}{k_{\hat\nu}^{1/2}}\,f_{\mu}=
\IM^{1/2}fk^{-1/2}
\eea
hence
\beq \label{eq:rR}
r_R=\IM^{-1/2}\,f^{-1}k^{1/2}g_R^{-1}k_R^{-1/2}f_R\ .
\eeq
Using (\ref{S}) we obtain $r_L=S(r_R)$.
 
In order to verify $\pi_M(\hat e)=E^L_M$ one uses
$\pi_M(\varphi)y=(\varphi\la ys^{1/2})s^{-1/2}$, $\varphi\in\duA,\
y\in A$, and the fact that within $A\cros\duA$ one has the
identification $\hat q_L=r_R$. Thus
\[
\pi_M(\hat e)y=r_R^{1/2}E^L(r_R^{1/2}ys^{1/2})s^{-1/2}=
E^L(yr_R^{1/2}s^{1/2})r_R^{1/2}s^{-1/2}=E^L(yr_R)
\]
where we utilized $\theta_{E^L}(r_R)=r_R$ and that $r_Rs^{-1}\in
A^L$. The proof of $\pi'_M(\hat e)=E^R_M$ goes analogously using
$\pi'_M(\varphi)y=(ys^{1/2}\ra \duS^{-1}(\varphi))s^{-1/2}$ and
the fact that within $\duA\cros A$ we have $\hat q_R=r_L$.
It is now a standard consequence \cite{G-H-J} that the relations
(\ref{TLJ1}) and (\ref{TLJ3}) hold true, using also the duality
principle for those involving $e$. The Temperley-Lieb relations
(\ref{TLJ2}-b) in turn follow from the fact that $e\in A$ and
$\hat e\in\duA$ after the reader have checked
that $E^L_M(e)=\IM^{-1}=E^R_M(e)$. \qed
 
Being the Jones extension $A\cros\duA=A\hat e A$, in particular
for
every $\varphi\in\duA\subset A\cros\duA$ there exist $a_i,b_i\in
A$ such that $\varphi=\sum_i a_i\hat e b_i$. In order to obtain a
concrete expression use the quasibasis $\sum_i u_i\o v_i=
S(h\c)\o\IM q_R h\cc$ of $E^L_M$:
\bea
\1&=&\sum_i u_i\hat e v_i\ =\ u_i\hat h r_R v_i\nn
\varphi=\varphi\1&=&\sum_i (\varphi\c\la u_i)\varphi\cc\hat h r_R
v_i=\sum_i (\varphi\c\la u_i)\duPL(\varphi\cc)\hat h r_Rv_i=\nn
&=&\sum_i(\varphi\la u_i)\hat h r_R v_i
=(\varphi\la S(h\c))r_R^{-1/2}\hat
er_R^{-1/2}g_R^{-2}h\cc             \label{fi}
\eea
This will be used to prove the following
\begin{coro}
Let $\tau_M^W\colon A\cros\duA\to\C$ be the trace associated to
$\tau_M\colon A\to\C$ by the basic construction for $A^L\subset
A$. Similarly, let $\tau_M^{W'}\colon A\cros\duA$ be the trace
associated to the Markov trace $t_M$ of $\duA$ by the basic
construction for $\duA^R\subset \duA$. Then
$\tau_M^W=\tau_M^{W'}$.
\end{coro}
Notice that, as a consequence of this, the restriction to
$\duA$ of the conditional expectation $\hat E^L_M\colon
A\cros\duA\to A$ defined by $\hat E^L_M(x\hat e y)=x\IM^{-1}y$
coincides with the Markov conditional expectation previously
denoted by $\hat E^L_M$. Together with the analogue statement for
$E^R_M\colon A\cros \duA\to \duA$, this means commutativity of the
Markov conditional expectations around the squares of Figure 2.
 
\Proof By definition $\tau_M^W(x\hat e
y)=\IM^{-1}\bra\tau_M,xy\ket$ for $x,y\in
A$ and thus, by cyclicity of the trace, it is sufficient to prove
that $\tau_M^W(e\varphi)=\IM^{-1}\bra\varphi,t_M\ket$ for
$\varphi\in\duA$.
\beanon
\tau_M^W(e\varphi)&=&\sum_i\tau_M^W(e(\varphi\la u_i)r_R^{-1/2}
\hat e r_R^{1/2}v_i)=\IM^{-1}\sum_i\tau_M(e(\varphi\la u_i)v_i)=\\
&=&\bra\tau_M,e(\varphi\la S(h\c))q_Rh\cc\ket=\bra\tau_M,hq_R^{1/2}
(\varphi\la S(h\c))q_Rh\cc q_R^{1/2}\ket
\eeanon
where we used the identity $r_Rq_R=\IM^{-1}g_R^{-2}$. Inserting
here the calculation
\[
\tau_M\ra h=h\la\tau_M=h\la \hat h\hat g_L^{-2}\hat k^{-1}\hat
f^2=\hat k^{-1}\hat f^2
\]
we obtain
\beanon
\tau_M^W(e\varphi)&=&\eps(q_R^{1/2}(\varphi\la S(h\c))q_Rh\cc
q_R^{1/2}k_R^{-1}f_R^2)=\bra\varphi,S(h\c)\ket\ \eps(q_R^{1/2}
\PR(q_Rh\cc)\\
&&q_R^{1/2}k_R^{-1}f_R^2)=
\bra\varphi,S(h\c)\ket\ \eps(q_Rh\cc
q_Rk_R^{-1}f_R^2)=\bra\varphi, S(q_Rhq_Rk_R^{-1}f_R^2)\ket=\\
&=&\bra\varphi,f_L^2k_L^{-1}q_Lhq_R\ket=\IM^{-1}\bra\varphi,t_M\ket
\eeanon
where in the last equation we took into account the dual of the
formula (\ref{RN of tau_M/Haar}).\qed
 
The restrictions of the Markov trace $\tau_M^W$ onto various
subalgebras of $A\cros\duA$ have trace vectors as listed below.
\[
\begin{array}{c|c|r}
\mbox{Subalgebra}&\mbox{minimal central projections}&\mbox{trace
vector}\\ \hline
Z& z_{\hat\nu}\,,\ \hat\nu\in\Vac\duA
&\IM_{[\hat\nu]} f_{\hat\nu}^2 \\
\hat Z& z_\mu\,,\ \mu\in\Vac A
&\IM_{[\mu]} f_{\mu}^2 \\
A^L& e_a^L\,,\ a\in \Sec A^L
&\IM_{[a]}^{1/2} f_{a^L} d_a f_{a^R} \\
A^R\equiv \duA^L& e_b^R\,,\  b\in \Sec A^R
&\IM_{[b]}^{1/2} f_{b^L} d_b f_{b^R} \\
\duA^R&\hat e_a^R\,,\ a\in \Sec A^L
&\IM_{[a]}^{1/2} f_{a^L} d_a f_{a^R} \\
A&e_q\,,\ q\in \Sec A
&f_{q^L} d_q f_{q^R} \\
\duA&\hat e_{\hat q}\,,\ \hat q\in \Sec\duA
&f_{\hat q^L} d_{\hat q} f_{\hat q^R} \\
A\cros\duA&e_a^W\,,\ a\in\Sec A^L
&\IM_{[a]}^{-1/2} f_{a^L} d_a f_{a^R}
\end{array}
\]
Some comments on the idempotents $e_a^W$ are in order. From the
general theory of Jones extensions we know that $\Center
A\cros\duA\cong\Center A^L$, however for an unambiguous labelling
of the minimal central idempotents of $A\cros\duA$ with $a\in\Sec
A^L$ we need the "shift isomorphism"
\beq  \label{shift}
\Center A^L\ \ni\ z^L\quad\mapsto\quad
\sum_i\,u_iz^L\hat e v_i\ \in\ \Center A\cros\duA\ .
\eeq
Thus our definition is this
\beq
e_a^W\ :=\ \sum_i u_i e_a^L\hat e v_i=S(h\c)e_a^L
\hat h g_R^{-2}h\cc\ .
\eeq
It is important to remark that one would have obtained the same
result for $e_a^W$ using the canonical isomorphism $A^L\to\duA^R$,
$e_a^L\mapsto \hat e_a^R=e_a^L\la\du1$ and after that the other
shift isomorphism $\Center\duA^R\to\Center A\cros\duA$ associated
to the basic construction for $\duA^R\subset \duA$. As a matter of
fact the reader may check the following equality in $A\cros\duA$
valid for all $z^L\in\Center A^L$,
\beq
S(h\c)z^L\hat hg_R^{-2}h\cc\ =\
\hat h\c\hat g_L^{-2}h(z^L\la\du1)\duS(\hat h\cc)\ ,
\eeq
which expresses commutativity of the triangle consisting of the
two shift isomorphisms and of the canonical isomorphism
$\kappa_A^L$ of Lemma I.2.6.
Similarly, there is an unambiguous labelling of the minimal
central projections $e_b^{W'}$ of $\duA\cros A$ by the sectors $b$
of $A^R$.
 
The content of the next Proposition can be phrased as Frobenius
reciprocity in the underlying 2-category of the weak Hopf algebra.
\begin{prop}
For $a,a'\in \Sec A^L$, $b,b'\in\Sec A^R$, and $q\in\Sec A$ let
$N^{qa'}_a$ be the multiplicity of the simple algebra $e_qA\o
\hat e_{a'}^R\duA^R$ in the simple algebra $e_a^W(A\cros\duA)$ and
let $N^{b'q}_b$ be the multiplicity of the simple algebra
$\hat e_{b'}^L\duA^L
\o e_q A$ in $e_b^{W'}(\duA\cros A)$. As before, $N_q^{ab}$
denotes the inclusion matrix of $A^L\amalgo{Z}A^R\subset A$. Then
\beq
N_a^{q\bar b}\ =\ N_q^{ab}\ =\ N_b^{\bar aq}
\eeq
where $a\mapsto\bar a$ and $b\mapsto \bar b$ are the mutually
inverse bijections induced by the antipode restricted to $\Center
A^{L/R}$, respectively.
\end{prop}
\Proof We will content ourselves with proving the first equality.
The subalgebra in $A\cros\duA$ generated by $A$ and $\duA^R$ is
the amalgamated tensor product $A\amalgo{\hat Z}\duA^R$ with
minimal central projections $e_q\hat e_a^R$ where $q^R=a^L$. The
inclusion matrix of $A\amalgo{\hat Z}\duA^R\subset A\cros\duA$ can
therefore be computed as follows.
\beanon
N_a^{q\bar b}&=&\frac{\IM^{1/2}_{[a]}}{t_a}\,\tau_M^W(e_a^W\,
e_q\hat e_{\bar b}^R)\frac{1}{n_qn_{\bar b}}=
\frac{\IM^{1/2}_{[a]}}{n_qt_an_b}\,\tau_M^W(S(h\c)e_a^L\,\hat e
\,e_b^Rr_R^{-1}g_R^{-2}h\cc e_q)=\\
&=&\frac{\IM^{-1/2}_{[a]}}{n_qt_an_b}\,\tau_M(S(h\c)e_a^L
e_b^Rr_R^{-1}g_R^{-2}h\cc e_q)=\\
&=&\frac{\IM^{-1/2}_{[a]}}{n_qt_an_b}\,\frac{1}{\tau_q}\tau_M(e_q)
\tr_q(g^{-1}e_a^Le_b^Rr_R^{-1}g_R^{-2})=
\frac{f_{q^L}k_{q^L}^{-1/2}f_{b^L}k_{b^L}^{-1/2}}{t_an_b}
\tr_q(g_L^{-1}e_a^Le_b^R)=\\
&=&\frac{k_{a^L}^{-1/2}k_{a^R}^{-1/2}}{d_an_b}\,N_q^{ab}\tr_a(g_L^{-1})
\tr_b(\1)\ =\ N_q^{ab}
\eeanon
\qed
\begin{coro}
The restriction of an irreducible representation $D_q$ of $A$ onto
the subalgebra $z_\mu^LA^L$ is either the zero representation (if
$q^L\neq \mu$) or a faithful representation (if $q^L=\mu$). Thus
the inclusion matrix $\Lambda$ of $A^L\subset A$ satisfies
\beq
\Lambda_{aq}\ >\ 0\quad\Leftrightarrow\quad a^L\ =\ q^L\ .
\eeq
\end{coro}
\Proof
If $a^L=q^L$ then $\hat e_{\bar a}^Le_q$ is a non-zero projection
in $\duA\cros A$ (due to the intersection $\duA^L\cap A=\hat Z$).
Hence there exists a $b$ such
that $N^{\bar aq}_b>0$. It follows from Frobenius reciprocity that
$N_q^{ab}>0$, i.e. $\Lambda_{aq}=\sum_b N_q^{ab} > 0$. \qed
 
For a pure $C^*$-WHA this means that every representation
represents $A^L$ (and $A^R$) faithfully. Even in the non-pure case
the maximal possible faithfulness is attained which is still
compatible with the groupoidlike sector composition.
 
\subsection{Pairing formula}
To conclude the general analysis we return to the beginnings and
give an expression of the canonical pairing in terms of the Markov
trace and of the Temperley-Lieb-Jones projections. This formula
can be the starting point of the reconstruction of a WHA from a
given inclusion data.
 
\begin{thm}
With $\tau_M$ denoting the Markov trace on the Weyl algebra
$A\cros\duA$ the canonical pairing of $\varphi\in\duA$ and $x\in A$
can be written as
\beq
\bra\varphi,x\ket\ =\ \tau_M(\IM^{3/2}\cdot \hat e\,e\,\varphi
\,{g''_L}^{1/2}x\,{g''_R}^{1/2}\,)
\eeq
where we introduced the notation
\beabc
g''_L&:=&f_L^{-1}g'_Lf^{-1}\ =\
f_L^{-1}k_L^{1/2}g_Lk^{1/2}f^{-1}\\
g''_R&:=&f^{-1}g'_Rf_R^{-1}\ =\
f^{-1}k^{1/2}g_Rk_R^{1/2}f_R^{-1}
\eeabc
\end{thm}
\Proof At first compute $E^L_M(ex)=E^L(r_Rq_L^{1/2}hq_L^{1/2}x)=
q_L^{1/2}h\c q_L^{1/2}x\c \bra\hat h,r_Rh\cc x\cc\ket=
q_L^{1/2}r_Lh\c S^{-1}(x\cc)q_L^{1/2}x\c \bra\hat h,h\cc\ket=
\IM^{-1}q_L^{-1/2}S^{-1}(x\cc)q_L^{1/2}x\c$ and then apply this
together with (\ref{fi}) to obtain
\beanon
\hat E_M^L(\hat ee\varphi)
&=&\hat E_M^L(\hat ee(\varphi\la S(h\c))r_R^{-1/2}\hat e
  r_R^{-1/2}g_R^{-2}h\cc)=\\
&=&\hat E_M^L(\hat eE^L_M(e(\varphi\la S(h\c))r_R^{-1/2})
  r_R^{-1/2}g_R^{-2}h\cc)=\\
&=&\IM^{-2}q_L^{-1/2}r_L^{-1/2}h\cc q_L^{1/2}S(h\ccc)r_R^{-1/2}
  g_R^{-2}h_{(4)}\bra\varphi,S(h\c)\ket=\\
&=&\IM^{-2}q_L^{-1/2}r_L^{-1/2}r_R^{-1/2}g_R^{-2}(h\ra\duS(\varphi))
  q_L^{1/2}=\\
&=&\IM^{-1}gq_R^{1/2}(h\ra\duS(\varphi))q_L^{1/2}\ .
\eeanon
Therefore
\beq       \label{pre-pairing}
\tau_M(\hat ee\varphi x)=\bra\tau_M,\hat E^L_M(\hat
ee\varphi)x\ket=\bra\tau_M,\IM^{-1}(h\ra\duS(\varphi))q_L^{1/2}x
q_R^{1/2}g\ket\ .
\eeq
The next step is to express the canonical pairing via the
canonical trace $\tau$,
\beanon
\bra\varphi,x\ket&=&\bra\duS(\varphi),(\hat h\la h)g_L^{-2}S^{-1}
(x)\ket=\\
&=&\bra\duS(\varphi),h\c\ket\,\bra\hat h,\,h\cc xg_R^{-2}\ket=\\
&=&\bra\tau,\,(h\ra\duS(\varphi))xg\ket\ .
\eeanon
Taking into account the relation (\ref{RN of tau_M/Haar}) of
$\tau_M$ to $\tau$ and then (\ref{pre-pairing}) we obtain
\beanon
\bra\varphi,x\ket&=&\bra\tau_M,(h\ra\duS(\varphi))xgf_L^{-1}
k_L^{1/2}k_R^{1/2}f_R^{-1}\ket=\\
&=&\tau_M(\hat ee\varphi\,\IM f_L^{-1}k_L^{1/2}q_L^{-1/2}x
q_R^{-1/2}k_R^{1/2}f_R^{-1})=\\
&=&\tau_M(\IM^{3/2}\hat ee\varphi
(g''_L)^{1/2}x(g''_R)^{1/2}\,)
\eeanon
\qed

\section{Special cases}
\subsection{Weak Kac algebras}
 
Weak Kac algebras (WKA) or, what is the same, generalized Kac
algebras
of \cite{Yama} are precisely the weak $C^*$-Hopf algebras that
have involutive antipodes: $S^2=\id$ \cite{NV1}. If $A$ is a WKA
then its dual weak $C^*$-Hopf algebra $\duA$ is also a WKA.
\begin{lem}  \label{lem WKA}
The following conditions for a $C^*$-WHA $A$ are equivalent.
\begin{description}
\item[i)] $A$ is a WKA, i.e. $S^2=\id$,
\item[ii)] $\hat h$ is a trace on $A$,
\item[iii)] $\eps(g_R^{-2})=\dim A$.
\end{description}
\end{lem}
\Proof Equivalence of {\bf (i)} and {\bf (ii)} has already been
proven in \cite{NV1}. For completeness we give here an independent
argument. At first we recall Subsection I.4.3 that the Haar
measure is $\hat h=g_Lg_R\la \tau$ where the trace $\tau$ has
trace vector $\tau_q=\tr_q g\equiv\tr_q g^{-1}$.
 
{\bf (i) $\Rightarrow$ (ii)} $S^2=\Ad_g=\id$ implies $g=\1$ by
uniqueness of the canonical grouplike element (Proposition I.4.4).
Then $\hat g=\du1$ and therefore $\theta_{\hat h}(x)=\hat g\la
x\ra\hat g=x$ by (I.4.29), i.e. $\hat h$ is a trace.
 
{\bf (ii) $\Rightarrow$ (i)} Using the tracial property of $\hat
h$
and Proposition I.4.9 we have $\hat h\c\o\hat h\cc=\hat h\cc\o\hat
h\c=\hat h\c\o\hat g\hat h\cc\hat g$. By nondegeneracy of $\hat h$
this implies $\varphi=\hat g\varphi\hat g$ for all
$\varphi\in\duA$, in particular $\hat g^2=\du1$. Since $\hat g\geq
0$, we obtain $\hat g=\du1$ and $\duS^2=\id$. Taking transpose,
$S^2=\id$ follows.
 
{\bf (i) $\Leftrightarrow$ (iii)} For $x^L\in A^L$ we have
$\eps(x^L)=\bra\hat h,x^L\ket=\tau(x^Lg_Lg_R)$ implying two
interesting identities,
\bea
\eps(g_R^{-1})&=&\tau(g_R)\\
\eps(g_R^{-2})&=&\sum_q\tau_q^2
\eea
The first one is useful in examples to determine $g_R$ once $g$
is known. The second one together with the inequality
\beq
\tau_q^2=(\tr_q\, g)(\tr_q\,g^{-1})\geq n_q^2   \label{tau>n}
\eeq
shows that $\eps(g_R^{-2})\geq\dim A$ and equality holds iff
$\tau_q=n_q$, $q\in\Sec A$, which in turn is equivalent to that
$g\in\CA$ by the well known property of the inequality
(\ref{tau>n}). Clearly $g\in\CA$ iff $S^2=\id$.\qed
 
The main result of this subsection is the following
\begin{thm}
Let $A$ be a weak Kac algebra. Then the indices of the Markov
and of the Haar conditional expectations coincide and take an
integer value on each hypercentral block,
\beq
\IH_H\ =\ \IM_H\ \in\ \N\,,\quad H\in\Hyp A\ .
\eeq
Moreover, for each vacuum $\mu\in\Vac A$ the dimension of $z^L_\mu
A$ is divisible by that of $z^L_\mu A^L$ and their ratio is
$\IH_{[\mu]}$, hence constant over the hypercentral block.
\end{thm}
 
\Proof By Lemma \ref{lem WKA} ii) the Haar state is tracial.
Therefore $E^L$ is a trace preserving conditional expectation the
index of which belongs to the common centers of $A$ and $A^L$. Now
Scholium A.4 implies that $E^L$ is the Markov conditional
expectation $E^L_M$, hence $\IH=\IM$. It remains to show that this
common index is
\beq      \label{deltaKac}
\IM_{[\mu]}\ =\ {\dim z_\mu^L A\over\dim z^L_\mu A^L}
\eeq
and then Lemma A.5 of the Appendix will imply that
$I_H=\IM_H$, $H\in\Hyp A$
are integers. For that purpose, and also for mere curiosity, we
compute the quantities $k_\mu$, $d_q$, $g_L$, $f_\mu$, and
$d_{\mu\nu}$ for WKA's.
Since $\1\cc S(\1\c)=\1$, $\eps|_{A^L}=\tr^L$, the
left regular trace of $A^L$ by Lemma \ref{RNder}.
Furthermore $\tau_q=n_q$ by (\ref{tau>n}). Thus we have
\bea
k_\nu&=&\eps(z^L_\nu)=\sum_{a\in\Sec A^L,a^L=\nu}\ (n^L_a)^2\ =
\ \dim(z^L_\nu A^L)\\
d_q&=&{\tau_q\over(k_{q^L}k_{q^R})^{1/2}}\ =
\ {n_q\over(k_{q^L}k_{q^R})^{1/2}}\ .
\eea
The modular automorphism of the Haar functional on $A$ is the
identity therefore $g_Lg_R\in\Center A$ by Proposition I.4.14 i).
But $g=\1$ implies $g_L=g_R$ therefore
$g_L\in A^L\cap A^R\cap\Center A=\Hypercenter A$.
Now Proposition 4.3 immediately
gives the Haar index
\beq
I=g_L^{-2}\ .
\eeq
The traces $\hat h$ and $\tau_M$ having the same trace preserving
conditional expectations onto $A^L$ and onto $A^R$, too,
may differ only in a hypercentral Radon-Nikodym
derivative. Comparing this to Eqn (\ref{RN of tau_M/Haar}) we see
that $f_Lk_L^{-1/2}k_R^{-1/2}f_R\in\Hypercenter A$
which is possible only if $f_Lk_L^{-1/2}$ is itself hypercentral,
due to the fullness of the hypercentral blocks. Taking into
account the normalization (\ref{normalize tau_M}) this
hypercentral element can be determined and yields the expression
\beq  \label{fKac}
f_\mu\ =\ \left({k_\mu\over\dim z_{[\mu]} A}\right)^{1/2}\,,\quad
\mu\in \Vac A\ .
\eeq
Therefore the Markov trace has trace vector
\beq
t_q\ =\ {n_q\over\dim z_H A}\ ,\qquad q\in H\,,\ H\in\Hyp A\ .
\eeq
This means that $\tau_M$ restricts to the normalized regular trace
on each hypercentral block $z_HA$.
The regular dimension matrix
\beq
d_{\mu\nu}\ =\ k_\mu^{-1/2}k_\nu^{-1/2}\sumq \ n_q^2
\eeq
has $f$ as its Perron-Frobenius eigenvector. Inserting \ref{fKac}
into the eigenvalue equation one obtains
\beq
\sum_{q\in\Sec A,\ q^L=\mu}\ n_q^2\ =\ k_\mu\cdot\IM_{[\mu]}
\eeq
which proves (\ref{deltaKac}) and the Theorem.
Especially for pure weak Kac algebras
we obtain that $\dim A$ is divisible by $\dim A^L$. \qed

\subsection{The $C^*$-WHA $B\o B^{op}$}
 
In the Appendix of I. we have shown that any separable algebra
$B$ together with a nondegenerate functional $E$ of
index $1$ determines a WHA structure on $B\o B^{op}$. We develop
further this construction in case when $B$ is a finite dimensional
$C^*$-algebra and compute the quantities introduced in this paper.
 
Let $B\cong \oplus_\mu M_{n_\mu}$ with a set of matrix units
$e_\mu^{ij}$ and minimal central projections $e_\mu$. We define
the trace $\tr$ on $B$ by setting $\tr e_\mu=n_\mu$ and will
also use the traces $\tr_\mu(x):=\tr(e_\mu x)$.
The nondegenerate functional $E(x)=\tr(\gamma^2 x)$ is given in
terms of a positive invertible $\gamma\in B$ satisfying
$\tr_\mu(\gamma^{-2})=1$ for all $\mu\in\Sec B$. The structure
maps of the $C^*$-WHA $A=B\o B^{op}$ are the following (cf I.
Appendix):
\bea
\cop(x\o y)&=&\sum_\mu\sum_{ij}\ (x\o e_\mu^{ij}\gamma^{-1})\
\o\ (\gamma^{-1}e_\mu^{ji} \o y)\ ,\\
\eps(x\o y)&=& \tr(\gamma^2 xy)\ ,\\
S(x\o y)&=& y\o\gamma^2 x\gamma^{-2}\ .
\eea
The left and right subalgebras are $A^L=B\o\1$, $A^R=\1\o B$.
The sectors of $A$ are pairs $(\mu,\nu)$ of sectors of $B$. All
$(\mu,\nu)$ is either a vacuum sector (if $\mu=\nu$) or a soliton
sector (if $\mu\neq\nu$). The dual $\duA$ is a simple algebra with
a single sector denoted $\circ$.
 
Using the Definitions I.3.1 and I.3.24 the reader may check that
the element
\beq
h\ =\ \sum_\mu\ {1\over \Gamma_\mu}\ \sum_{ij}\ e_\mu^{ij}\gamma\o
    \gamma e_\mu^{ji}
\eeq
is the Haar integral in $A$, where
$\Gamma_\mu:=\tr_\mu(\gamma^2)$. Introducing the notation
$\Gamma:=\sum_\mu \Gamma_\mu e_\mu$ the canonical grouplike
element can be written as
\beq\label{BoB g}
g\ =\ \Gamma^{-1/2}\gamma^2\ \o\ \gamma^{-2}\Gamma^{1/2}\ .
\eeq
Hence the trace vector of the canonical trace $\tau$ is
$\tau_{(\mu,\nu)}=(\tr_\mu\o\tr_\nu)(g)=\sqrt{\Gamma_\mu\Gamma_\nu}$.
Using the identity $\eps(g_R^{-1})=\tau(g_R)$ (\ref{BoB g})
implies
\beabc
g_L&=&\frac{1}{(\sum_\mu\Gamma_\mu)^{1/2}}\Gamma^{-1/2}\gamma^2\o\1\\
g_R&=&\frac{1}{(\sum_\mu\Gamma_\mu)^{1/2}}\1\o\gamma^2\Gamma^{-1/2}
\eeabc
In this example $Z^L=\Center A^L$ has minimal idempotents
$z_\mu^L=e_\mu\o\1$, therefore the function $k_\mu=\Gamma_\mu$ and
$k_L=\Gamma\o\1$, $k_R=\1\o\Gamma$. We obtain for the standard
metric the expression $g'=\gamma^2\o\gamma^{-2}$ and the
dimensions of all of the sectors are $d_{(\mu,\nu)}=1$. Taking
into account that $k_{\circ}=\eps(\1)=\sum_\nu\Gamma_\nu$ the
dimensions of the sectors of $A^L$ and $A^R$ are also trivial:
$d_\mu=1$. Hence the left regular dimension matrices of $A^L$,
$A$, $A^R$, and $\duA$ are
\beq
d_{\mu\circ}=n_\mu\,,\quad d_{\mu\nu}=n_\mu n_\nu\,,\quad
d_{\circ\nu}=n_\nu\,,\quad d_{\circ\circ}=\sum_\nu n_\nu^2\,,
\eeq
respectively. Hence the Perron-Frobenius eigenvectors are
$f_\mu=n_\mu/\dim B$ and $f_\circ=1/\sqrt{\dim B}$.
The common eigenvalue, which is the Markov index of the inclusions
$A^{L/R}\subset A$, is $\IM=\sum_\nu n_\nu^2=\dim B=\sqrt{\dim
A}$. For generic choices of $\gamma$ the three traces
$\tau$, $\tau_S$, and $\tau_M$ are different:
\bea
\tau(x\o y)&=&\sum_{\mu,\nu}\ \sqrt{\Gamma_\mu\Gamma_\nu}\
(\tr_\mu\, x)(\tr_\nu\, y)\\
\tau_S(x\o y)&=&(\tr\, x)(\tr\, y)\\
\tau_M(x\o y)&=&\sum_{\mu,\nu}\ \frac{n_\mu n_\nu}{(\dim B)^2}\
(\tr_\mu\, x)(\tr_\nu\, y)\ .
\eea
The Markov trace coincides with the (normalized) trace in the left
regular representation of $A$. The Haar functional $\hat h$ and
the Haar conditional expectation $E^L(x)=\hat h\la x$ are now easy
to evaluate,
\bea
\bra\hat h,x\o y\ket&=&\frac{(\tr\gamma^2 x)(\tr\gamma^2 y)}
{\tr\gamma^2}\ ,\\
E^L(x\o y)&=&x\o\1\frac{\tr \gamma^2 y}{\tr\gamma^2}\ .
\eea
For the Haar index $I=\Ind E^L$ one obtains the scalar
$I=\1_A\tr\gamma^2$.

\appendix
\sec{On the index of finite dimensional inclusions}
 
The following results may belong to the standard part of the
theory
of inclusions of multimatrix algebras \cite{G-H-J}, although it is
difficult to find them in the form presented here, mainly because
we have been using Watatani's ring theoretical notion of
index \cite{Watatani}.
 
\begin{scho}
Let $A\subset B$ be a unital inclusion of finite dimensional
$C^*$-algebras and let $\varphi\colon B\to \C$ be a faithful
positive linear functional.
Define the $A$-module maps $E\colon\, _AB\to\,_AA$ and
$F\colon B_A\to A_A$ respectively by the formulae
$$
\varphi(aE(b))=\varphi(ab)\,,\quad
\varphi(F(b)a)=\varphi(ba)\qquad a\in A,b\in B\ .
$$
Then the following statements are equivalent:
\begin{description}
\item[i)] $\theta(A)\subset A$
\item[ii)] $E\circ\theta=\theta\circ F$
\item[iii)] $E=F$
\item[iv)] $E$ is a conditional expectation.
\end{description}
\end{scho}

If the above equivalent conditions hold then $E_\varphi=E$ will be
called the {\em $\varphi$-preserving conditional expectation}.
 
Now let $\varphi$ and $\psi$ be faithful positive functionals on
$B$ such that $\theta_\varphi(A)=A$ and $\theta_\psi(A)=A$. We
define the (left) Radon-Nikodym derivatives of $\varphi$ w.r.t
$\psi$ and of $E_\varphi$ w.r.t $E_\psi$, respectively by
\beq
\varphi(b)=\psi(s b)\,,\quad E_\varphi(b)=E_\psi(rb)\,,
\quad b\in B\ .
\eeq
$r$ can be computed from $s$ by the following formulae
\beq                                      \label{r(s)}
r\ =\ E_\varphi(s^{-1})s\ =\ E_\psi(s)^{-1}s\ =\ s
E_\psi(s)^{-1}\ .
\eeq
As a matter of fact $\varphi(E_\varphi(b)a)=\varphi(ba)=
\psi(s ba)=\psi(E_\psi(s b)a)=\varphi(E_\varphi(s^{-1})
E_\psi(s b)a)$, implying the first equality in (\ref{r(s)}).
The second follows from the first because $E_\psi(r)=\1$ and the
third follows from the second since $r\in A'\cap B$
\cite{Watatani}. One can also show easily that the modular
automorphisms are related by
\beq
\theta_\varphi=\theta_\psi\circ\Ad_s\,,\quad
\theta_{E_\varphi}=\theta_\varphi|_{A'\cap B}=
\theta_{E_\psi}\circ\Ad_r|_{A'\cap B}\ ,
\eeq
and that $r$ and $s$ commute. If $\psi=\tau$ is tracial then
$r$ is positive and we have
\beq    \label{E_fi(E_tau)}
E_\varphi(b)=E_\tau(rb)=E_\tau(br)=E_\tau(r^{1/2}br^{1/2})\,,
\quad b\in B\ .
\eeq
For a fixed faithful trace $\tau$ and for arbitrary $\varphi$ and
$\psi$ as above let $s_\varphi,\ s_\psi$ and $r_\varphi,\
r_\psi$ be the corresponding Radon-Nikodym derivatives w.r.t
$\tau$ and $E_\tau$, respectively. Then one has the following
manifestly positive expressions
\bea
\varphi(b)&=&\psi(s_\psi^{-1/2}s_\varphi^{1/2}\,b\,s_\varphi^{1/2}
s_\psi^{-1/2})\ ,\\
E_\varphi(b)&=&E_\psi(r_\psi^{-1/2}r_\varphi^{1/2}\,b\,r_\varphi^{1/2}
r_\psi^{-1/2})\ .
\eea
 
In order to study the index of various conditional expectations
we need the inclusion data $A\subset B$ explicitely:
$A\cong\oplus_\alpha M_{n_\alpha}$,
$B\cong\oplus_\beta M_{m_\beta}$, and inclusion matrix
$\Lambda=[\Lambda_{\beta\alpha}]$. Then there exists a set
$\{e_\beta^{IJ}\,|\, I,J\in\I_\beta,\ \beta\in\Sec B\,\}$ of
matrix units for $B$ where the index set $\I_\beta$ consists of
triples $I=(a,\alpha,i)$ where $\alpha\in\Sec A$, is such that
$\Lambda_{\beta\alpha}>0$, $i=1,\dots,\Lambda_{\beta\alpha}$, and
$a=1,\dots, n_\alpha$. An arbitrary conditional expectation
\beq
E\colon B\to A\ ,\quad E(e_\beta^{a'\alpha'i',i\alpha a})=
\delta^{\alpha'\alpha}\ \Phi^{i'i}_{\beta\alpha}\ e_\alpha^{a'a}
\eeq
can be uniquely characterized by positive elements
$\Phi_{\beta\alpha}\in M_{\Lambda_{\beta\alpha}}$ satisfying
$\sum_\beta \tr\Phi_{\beta\alpha}=1\,,\ \forall \alpha$. Here
\beq
e_\alpha^{a'a}\ =\ \sum_\beta\sum_k e_\beta^{a'\alpha k,k\alpha a}
\qquad \alpha\in\Sec A,\ a',a=1,\dots, n_\alpha
\eeq
are matrix units for $A$. A faithful conditional expectation $E$
corresponds to having $\Phi_{\beta\alpha}$ invertible whenever
$\Lambda_{\beta\alpha}>0$. In the latter case we choose invertible
$C_{\beta\alpha}$ such that
$\Phi_{\beta\alpha}=C_{\beta\alpha}C_{\beta\alpha}^*$. Then it is
straightforward to verify that the set of elements
\beq
b_\beta^{a'\alpha'i',i\alpha a}\ :=\ \sum_j\
e_\beta^{a'\alpha'i',j\alpha a}\,{1\over\sqrt{n_{\alpha}}}\,
(C_{\beta\alpha}^{-1*})^{ji}
\eeq
form a quasibasis of $E$, i.e. $\sum_\beta\sum_{I,J\in\I_{\beta}}
\,b_{\beta}^{IJ}E(b_\beta^{IJ*}x)=x,\ \forall x\in B$. Therefore
\beq   \label{Index E}
\Ind E\ =\ \sum_{I\beta J}\ b_\beta^{IJ}b_\beta^{IJ*}\ =\
\sum_{\beta\alpha}\, e_{\beta}\cdot \tr\Phi_{\beta\alpha}^{-1}
\eeq
where $e_\beta=\sum_Ie_\beta^{II}$. Let
$\{f_{\beta\alpha}^i\,|\,i=1,\dots,\Lambda_{\beta\alpha}\,\}$ be
the eigenvalues of $\Phi_{\beta\alpha}$. Then
\bea
\sum_\beta\sum_i\ f_{\beta\alpha}^i&=&1\\
\Ind E&=&\sum_\beta\,e_\beta\cdot\sum_\alpha\sum_i\,{1\over
f_{\beta\alpha}^i}\ .
\eea
Let us choose a set $\{w_{\alpha}\}$ of positive
numbers. Then the inequality between arithmetic and harmonic means
weighted by $\{w_{\alpha}\}$ yields
\beq
{\sum_{\alpha i}\,{1\over f_{\beta\alpha}^i}\over
 \sum_\alpha \,\Lambda_{\beta\alpha} w_{\alpha}}\ \geq\
{\sum_\alpha\,\Lambda_{\beta\alpha} w_{\alpha} \over
 \sum_\alpha\sum_i\,f_{\beta\alpha}^iw_{\alpha}^2}
 \eeq
implying the estimate
\beq
\Ind E\ \geq\ \sum_\beta\ e_\beta\cdot{\left(\sum_\alpha\,
\Lambda_{\beta\alpha}w_{\alpha}\right)^2\over
\sum_\alpha\sum_i\, f_{\beta\alpha}^iw_{\alpha}^2}
\eeq
valid for all sequences $\{w_{\alpha}\}$ of positive numbers.
Equality holds here iff there exist numbers $u_\beta$ such that
$f_{\beta\alpha}^i=u_\beta/w_{\alpha}$ for all $i$.
 
\begin{lem}
Let $E\colon B\to A$ be a faithful conditional expectation over
the connected inclusion $A\subset B$ with inclusion matrix
$\Lambda$. Then
\beq                               \label{index bound}
\|\Ind E\|\ \geq\ \|\Lambda\|^2
\eeq
where $\|\,.\,\|$ denotes $L^2$-operator norm.
\end{lem}
\Proof $\|\Ind E\|=\max_\beta\,\sum_{\alpha,i}{1\over
f_{\beta\alpha}^i}\ \geq\
\max_\beta\,\frac{(\sum_\alpha\Lambda_{\beta\alpha}w_\alpha)^2}
{\sum_{\alpha,i}f_{\beta\alpha}^iw_\alpha^2}$ for all choices of
positive numbers $\{w_\alpha\}$. Using the identity
\[
\sum_\beta\sum_{\alpha,i}\ f_{\beta\alpha}^iw_\alpha^2\ =\
\sum_\alpha\ w_\alpha^2
\]
we obtain
\beanon
\|\Ind E\|&\geq&\left\{\max_\beta{\left(\sum_\alpha\Lambda_{\beta\alpha}
w_\alpha\right)^2\over\sum_{\alpha,i}f_{\beta\alpha}^iw_\alpha^2}
\right\}\
\sum_\beta\sum_{\alpha,i}f_{\beta\alpha}^iw_\alpha^2\ \cdot\
{1\over\sum_\alpha w_\alpha^2}\\
&\geq&\sum_\beta\
{\left(\sum_\alpha\Lambda_{\beta\alpha}w_\alpha\right)^2\over
\sum_\alpha w_\alpha^2}\ =\ {\|\Lambda w\|^2\over \|w\|^2}
\eeanon
for all vectors $w$ with positive entries. Choosing $w$ to be the
Perron-Frobenius eigenvector of $\Lambda^t\Lambda$ we obtain the
desired result.\qed
 
Now we turn to the special case of trace preserving conditional
expectations. Let $\tau\colon B\to\C$ be a faithful trace with
trace vector $\{t_\beta\}$, i.e.
$\tau(e_\beta^{IJ})=\delta^{IJ}t_\beta$ and $t_\beta>0$. Then the
$\tau$-preserving conditional expectation is the
unique $E\colon B\to A$ satisfying $\tau(ab)=\tau(aE(b))$ for
$a\in A$, $b\in B$.
The $\Phi$ matrices of this $E$ can be computed to be
\beq
\Phi_{\beta\alpha}^{i'i}\ =\ \delta^{i'i}\,{t_\beta\over s_\alpha}
\eeq
where $s_\alpha=\sum_\beta \, t_\beta\Lambda_{\beta\alpha}$ is the
trace vector of $\tau|_A$.  Inserting this into the general
formula (\ref{Index E}) gives
\beq \label{trace preserving Index E}
\Ind E\ =\ \sum_\beta\
e_\beta\,\cdot\,{\sum_{\beta'}(\Lambda\Lambda^t)_{\beta\beta'}
t_{\beta'}\over t_\beta}\ .
\eeq
 
\begin{defi}                     \label{Markov cond exp}
Let $A\subset B$ be a connected inclusion with inclusion matrix
$\Lambda$. Then the trace $\tau_M\colon B\to \C$ is called the
Markov trace for $A\subset B$ if its trace vector $t$ is the
Perron-Frobenius eigenvector of $\Lambda\Lambda^t$.
To make it unique we require $\tau_M(\1)=1$, i.e. $\sum_\beta
m_\beta t_\beta=1$. The $\tau_M$-preserving conditional
expectation $E_M\colon B\to A$ is called the Markov conditional
expectation.
 
If $A\subset B$ is not connected then let $\{z_\nu\}$ be the set
of minimal idempotents in $\Center A\cap\Center B$ and let
$\tau_\nu$ be the Markov trace of $z_\nu A\subset z_\nu B$. Then
the trace preserving conditional expectation associated to any
trace $\tau$ on $B$ the restriction of which to
the connected components $z_\nu A\subset z_\nu B$ are nonzero
multiples of $\tau_\nu$ is a unique conditional expectation
$E_M\colon B\to A$, called the Markov conditional expectation.
\end{defi}
 
Now a quick look at formula (\ref{trace preserving Index E})
yields the following
 
\begin{scho}
If $E$ is a trace preserving conditional expectation then
\beq
\Ind E\in\Center A\cap\Center B\quad\Longleftrightarrow\quad
E=E_M
\eeq
and the index of $E_M$ saturates the bound (\ref{index bound}),
i.e. for connected inclusions $\Ind E_M=\1\cdot \|\Lambda\|^2$.
\end{scho}
 
Finally we study connected inclusions for which the Markov index
takes its na\"{\i}ve value $\dim B/\dim A$. Notice that since
$\Lambda n=m$, we have the general estimate
\beq
\|\Lambda\|^2\ \geq\ {\|\Lambda n\|^2\over \|n\|^2}\ =\
{\dim B\over \dim A}\ .
\eeq
\begin{lem}
Let $A\subset B$ be a connected inclusion such that the norm of
the inclusion matrix satisfies $\|\Lambda\|^2=\dim B/\dim A$. Then
the Markov trace $\tau_M$ is the left regular trace on $B$ and the
Markov index $\|\Lambda\|^2$ is an integer. Hence $\dim B$ is
divisible by $\dim A$.
\end{lem}
 
\Proof Since ${\|\Lambda n\|^2\over\|n\|^2}=\|\Lambda\|^2$, $n$ is
the Perron-Frobenius eigenvector of $\Lambda^t\Lambda$. But then
$m=\Lambda n$ is the Perron-Frobenius eigenvector of
$\Lambda\Lambda^t$, thus
\beq
t_\beta={m_\beta\over\dim B}\qquad s_\alpha={n_\alpha\over\dim A}
\eeq
and $\tau_M$ and $\tau_M|_A$ are the normalized regular traces on
$B$ and $A$, respectively. Now the dimension vectors $m$ and $n$
satisfy the equations
\beq
\Lambda n\ =\ m\ ,\qquad \Lambda^t m\ =\ n\cdot I
\eeq
where $I=\|\Lambda\|^2$, the Markov index. The 2nd equation
implies that $n_\alpha I$ are integers therefore if $l$ denotes
the greatest common divisor of $\{n_\alpha\,|\,\alpha\in\Sec A
\,\}$ then $lI\in\Z$. Now the 1st equation implies that each
$m_\beta$ is divisible by $l$, too, hence $m'={m\over l}$ and
$n'={n\over l}$ are also integer vectors and satisfy $\Lambda^t
m'=n'\cdot I$. Therefore $I$ is an integer.\qed


\begin{thebibliography}{FroGab}
\addcontentsline{toc}{section}{\protect\numberline{}{References}}
\renewcommand{\b}{\bibitem}
\renewcommand{\baselinestretch}{.3}
\small
\medskip
\b{Westbury} J. W. Barrett, B. W. Westbury, Spherical
categories, preprint (1993)
\b{B} G. B\"ohm, Weak Hopf algebras and their application to
spin models, {\em PhD-thesis}, Budapest, 1997
\b{BSz} G. B\"ohm, K. Szlach\'anyi, A coassociative
$C^*$-quantum group with nonintegral dimensions, {\em Lett. Math.
Phys.} {\bf 35}, 437--456 (1996)
\b{BNSz} G. B\"ohm, F. Nill, K. Szlach\'anyi, Weak Hopf
Algebras I:  Integral Theory and the $C^*$-structure,
{\em math.QA/9805116}, to appear in J. Algebra
\b{Deligne-Milne} P. Deligne, J.S. Milne, Tannakian
categories, {\em Lecture Notes in Mathematics} Vol.900,
pp.101--228, Springer 1982
\b{EV} M. Enock, J.-M. Vallin, Inclusions of von Neumann
Algebras, and quantum groupoids, Inst. de Math. de
Jussieu, preprint No.156, 1988
\b{FGSzV} J. Fuchs, A. Ganchev, K. Szlach\'anyi, P. Vecserny\'es,
$S_4$-symmetry of $6j$-symbols and Frobenius-Schur indicators
in rigid monoidal $C^*$-categories,
{\em J. Math. Phys.} {\bf 40}, 408--426 (1999)
\b{Fredenhagen} K. Fredenhagen, Generalizations of the theory
of superselection sectors, in {\em Algebraic Theory of
Superselection Sectors}, ed. D. Kastler, World Scientific, 1989
\b{G-H-J} F.M. Goodman, P. de la Harpe, V.F.R. Jones:
{\em Coxeter Graphs and Towers of Algebras}, Springer 1989
\b{Longo} R. Longo, A duality for Hopf algebras and for
subfactors. I., {\em Commun. Math. Phys.} {\bf 159}, 133 (1994)
\b{Longo-Roberts} R. Longo, J. E. Roberts, A theory of
dimension, {\em K-theory} {\bf 11}, 103--159 (1997)
\b{MacLane} S. Mac Lane: {\em Categories for the Working
Mathematician}, Graduate Text in Mathematics 5, Springer-Verlag
1971
\b{NV1} D. Nikshych, L. Vainerman, Algebraic versions of a
finite-dimensional quantum groupoid,
{\em math.QA/9808054}
\b{Nill} F. Nill, Axioms for Weak Bialgebras,
{\em math.QA/9805104}
\b{NSzW} F. Nill, K. Szlach\'anyi, H.-W. Wiesbrock, Weak Hopf
algebras and reducible Jones inclusions of depth 2,
{\em math.QA/9806130}
\b{NSz} F. Nill, K. Szlach\'anyi, Quantum chains of Hopf
algebras and order-disorder fields with quantum double
symmetry, {\em hep-th 9507 174};
Quantum chains of Hopf algebras with quantum double cosymmetry,
{\em Commun. Math. Phys.} {\bf 187}, 159-200 (1997)
\b{Popa} S. Popa, Classification of amenable subfactors of
type II, {\em Acta. Math.} {\bf 172}, 163--255 (1994)
\b{Sz} K. Szlach\'anyi, Weak Hopf Algebras, in {\em Operator
Algebras and Quantum Field Theory}, eds. S. Doplicher, R. Longo,
J.E. Roberts, L. Zsid\'o, International Press, 1996
\b{Szymanski} W. Szymanski, Finite index subfactors and Hopf
algebra crossed products, {\em Proc. Amer. Math. Soc.} {\bf 120},
519 (1994)
\b{Watatani} Y. Watatani, Index for $C^*$-subalgebras,
{\em Memoirs of the Amer. Math. Soc.}, No. 424, 1990
\b{Wenzl} H. Wenzl, On the structure of Brauer's centralizer
algebras, {\em Annals of Math.} {\bf 128}, 173-193 (1988)
\b{Yama}T. Yamanouchi, Duality for generalized Kac algebras
and a characterization of finite groupoid algebras,
{\em J. Algebra} {\bf 163}, 9-50 (1994)
\b{Yetter} D. N. Yetter, Framed tangles and a theorem of
Deligne on braided deformations of Tannakian categories,
{\em Contemp. Math.} {\bf 134}, 325-349 (1992)
 
\end{thebibliography}
\end{document}